\documentclass[11pt]{elsarticle}

\usepackage{amsmath, amsthm, amssymb, mathabx}
\usepackage{graphicx}
\usepackage{color}
\usepackage{fullpage}

\begin{document}

\begin{frontmatter}

\title{Partitioning strategies for the interaction of a fluid with a poroelastic material based on a Nitsche's coupling approach}

\author[MA]{M. Buka\v{c}}
\ead{martinab@pitt.edu}
\author[MA] {I. Yotov}
\ead{yotov@math.pitt.edu}
\author[ME]{R. Zakerzadeh}
\ead{raz25@pitt.edu}
\author[ME]{P. Zunino\corref{cor1}}
\ead{paz13@pitt.edu}
\cortext[cor1]{Authors are listed in alphabetical order. Corresponding author (+1 412 624 9774)}

\address[MA]{Department of Mathematics, University of Pittsburgh, Pittsburgh, PA 15260,  USA}
\address[ME]{Department of Mechanical Engineering \& Materials Science, University of Pittsburgh, Pittsburgh, PA 15261,  USA}

\begin{abstract}
We develop a computational model to study the interaction of a fluid with a poroelastic material. The coupling of Stokes and Biot equations represents a prototype problem for these phenomena, which feature multiple facets. On one hand, it shares common traits with fluid-structure interaction. On the other hand it resembles the Stokes-Darcy coupling. For these reasons, the numerical simulation of the Stokes-Biot coupled system is a challenging task. The need of large memory storage and the difficulty to characterize appropriate solvers and related preconditioners for the equations at hand are typical shortcomings of classical discretization methods applied to this problem, such as the finite element method for spatial discretization and finite differences for time stepping. The application of loosely coupled time advancing schemes mitigates these issues, because it allows to solve each equation of the system independently with respect to the others, at each time step. In this work, we develop and thoroughly analyze a loosely coupled scheme for Stokes-Biot equations. The scheme is based on Nitsche's method for enforcing interface conditions. Once the interface operators corresponding to the interface conditions have been defined, time lagging allows us to build up a loosely coupled scheme with good stability properties. The stability of the scheme is guaranteed provided that appropriate stabilization operators are introduced into the variational formulation of each subproblem. The error of the resulting method is also analyzed, showing that splitting the equations pollutes the optimal approximation properties of the underlying discretization schemes. In order to restore good approximation properties, while maintaining the computational efficiency of the loosely coupled approach, we consider the application of the loosely coupled scheme as a preconditioner for the monolithic approach. Both theoretical insight and numerical results confirm that this is a promising way to develop efficient solvers for the problem at hand. 
\end{abstract}


\end{frontmatter}

%\title{Partitioning strategies for the interaction of a fluid with a poroelastic material based on a Nitsche's coupling approach}
%\author{M. Bukac, I. Yotov, R. Zakerzadeh, P. Zunino}
%\date{\today}
%
%\maketitle 

%>>>>>>>>>>>>>>>>>>>>>>>>>>>>>>>>>>>>>>>>>>>>>>>>>>>>>>>>>>>>>>>>>>>>>>>>>>>>>>>>>>>>>>>>>>>>>>>>>>>>>>>>>>>>>>>>>>>>>>>>>>>>>>>>>

% definitions
\def\DD{\boldsymbol D}
\def\uu{\boldsymbol u}
\def\nn{\boldsymbol{n}}
\def\hnn{\hat{\boldsymbol{n}}}
\def\ttau{\boldsymbol{t}}
\def\httau{\hat{\boldsymbol{t}}}
\def\hGamma{\hat{\Gamma}}
\def\bilinAA{{\mathcal A}}
\def\dt{d_\tau}
\def\dtt{d_\tau}
\def\interf{\Gamma}
\def\aalpha{}

%fluid domain
\def\domf{\Omega_f}
\def\hdomf{\hat{\Omega}_f}
\def\velf{\boldsymbol{v}_h}
\def\evelf{\widehat{\boldsymbol{e}}_{v,h}}
\def\hvelf{\hat{\boldsymbol{v}}_h}
\def\bff{\boldsymbol{\varphi}_{f,h}}
\def\hbff{\hat{\boldsymbol{\varphi}}_{f,h}}
\def\strf{\boldsymbol{\sigma}_{f,h}}
\def\hstrf{\hat{\boldsymbol{\sigma}}_{f,h}}
\def\strfv{\boldsymbol{D}_{f,h}}
\def\strfp{\pf}
\def\pf{p_{f,h}}
\def\epf{\widehat{e}_{p,f,h}}
\def\hpf{\hat{p}_{f,h}}
\def\bpf{\psi_{f,h}}
\def\hbpf{\hat{\psi}_{f,h}}

% porous matrix fluid equations
\def\domp{\Omega_p}
\def\hdomp{\hat{\Omega}_p}
\def\velp{\boldsymbol{q}_h}
\def\evelp{\widehat{\boldsymbol{e}}_{q,h}}
\def\hvelp{\hat{\boldsymbol{q}}_h}
\def\bfp{{\boldsymbol{r}}_h}
\def\hbfp{\hat{\boldsymbol{r}}}
\def\pp{p_{p,h}}
\def\epp{\widehat{e}_{p,p,h}}
\def\hpp{\hat{p}_{p,h}}
\def\bpp{\psi_{p,h}}
\def\hbpp{\hat{\psi}_{p,h}}
\def\velr{\boldsymbol{w}_h}
\def\velz{\boldsymbol{u}_h}
\def\hvelr{\hat{\boldsymbol{w}}_h}

% porous matrix solid equations
\def\strp{\boldsymbol{\sigma}_{p,h}}
\def\hstrp{\hat{\boldsymbol{\sigma}}_{p,h}}
\def\dispp{\boldsymbol{U}_h}
\def\edispp{\widehat{\boldsymbol{e}}_{U,h}}
\def\dotdispp{\dot{\boldsymbol{U}}_h}
\def\edotdispp{\dot{\widehat{\boldsymbol{e}}}_{U,h}}
\def\hdispp{\hat{\boldsymbol{U}}_h}
\def\bsp{\boldsymbol{\varphi}_{p,h}}
\def\dotbsp{\dot{\boldsymbol{\varphi}}_{p,h}}
\def\hbsp{\hat{\boldsymbol{\varphi}}_{p,h}}
\def\strep{\boldsymbol{\sigma}_{E,h}}
\def\hstrep{\hat{\boldsymbol{\sigma}}_{E,h}}

% membrane
\def\dispm{\boldsymbol{\eta}_h}
\def\edispm{\widehat{\boldsymbol{e}}_{\eta,h}}
\def\dotdispm{\dot{\boldsymbol{\eta}}_h}
\def\edotdispm{\dot{\widehat{\boldsymbol{e}}}_{\eta,h}}
\def\hdispm{\hat{\boldsymbol{\eta}}_h}
\def\bm{\boldsymbol{\xi}_h}
\def\dotbm{\dot{\boldsymbol{\xi}_h}}
\def\hbm{\hat{\boldsymbol{\xi}}_h}

% fem spaces
\def\femvf{\mathbf{V}_h^f}
\def\femvp{\mathbf{V}_h^p}
\def\fempf{Q_h^f}
\def\fempp{Q_h^p}
\def\fems{\mathbf{X}_h^p}
\def\dotfems{\dot{\mathbf{X}}_h^p}
\def\hfems{\hat{\mathbf{X}}_h^p}
\def\femm{\mathbf{X}_h^m}
\def\dotfemm{\dot{\mathbf{X}}_h^m}
\def\hfemm{\hat{\mathbf{X}}_h^m}

% bilinear forms
\def\bilinasp{a_{s}}
\def\bilinbsp{b_{s}}
\def\bilinm{a_m}
\def\bilinaff{{a}_{f}}
\def\bilinbff{{b}_{f}}
\def\bilinafp{{a}_{p}}
\def\bilinbfp{{b}_{p}}
\def\bilincfp{{c}_{p}}
\def\stabffp{{s}_{f,p}}
\def\stabffv{{s}_{f,v}}
\def\stabffq{{s}_{f,q}}

% energy
\def\energyf{E_{f,h}}
\def\energys{E_{p,h}}
\def\energym{E_{m,h}}

% error
\def\errorf{\mathcal{E}_{f,h}}
\def\errors{\mathcal{E}_{p,h}}
\def\errorm{\mathcal{E}_{m,h}}

% symmetry
\def\sym{\varsigma}

% rhs
\def\fff{{\cal F}}
\def\TT{{\cal T}}

%stability analysis
\def\epsfa{\widehat{\epsilon}_f^\prime}
\def\epsfab{\widecheck{\epsilon}_f^\prime}
\def\epsfaa{\epsilon_f^{\prime\prime}}
\def\epsfaaa{\widecheck{\epsilon}_f^{\prime\prime\prime}}
\def\epsfiv{\widehat{\epsilon}_f^{\prime\prime\prime}}

\def\epspa{\epsilon_p^\prime}
\def\epspaa{\epsilon_p^{\prime\prime}}
\def\epspaaa{\epsilon_p^{\prime\prime\prime}}

% convergence analysis
\def\yy{{\boldsymbol y}_h}
\def\yysn{{\boldsymbol y}_h(t_n)}
\def\yyi{\widetilde{\boldsymbol y}}
\def\yye{\widehat{\boldsymbol y}}
\def\yyin{\widetilde{\boldsymbol y}_h^n}
\def\yyen{\widehat{\boldsymbol y}_h^n}
\def\ein{\widetilde{\boldsymbol e}_h^n}
\def\een{\widehat{\boldsymbol e}_h^n}
\def\eeN{\widehat{\boldsymbol e}_h^N}
\def\zzh{{\boldsymbol z}_h}
\def\zz{{\boldsymbol z}_h}
\def\aas{{\cal A}_h}
\def\aasi{\widetilde{\cal A}_h}
\def\aasis{\widetilde{\cal A}_h^s}
\def\aase{\widehat{\cal A}_h}
\def\rsi{\widetilde{\cal R}_h}
\def\rse{\widehat{\cal R}_h}
\def\stab{{\cal S}_h}
\def\spliterr{{\cal R}_h}

\def\maas{{\mathrm{A}}_h}
\def\maasi{\widetilde{\mathrm{A}}_h}
\def\maasis{\widetilde{\mathrm{A}}_h^s}
\def\maase{\widehat{\mathrm{A}}_h}
\def\YY{{\boldsymbol Y}}
\def\ZZ{{\boldsymbol Z}}
\def\XX{{\boldsymbol X}}
\def\FF{{\boldsymbol F}}
\def\RR{{\boldsymbol R}}

\def\sumtau{\tau \sum_{n=1}^N}

% commands
\newcommand{\triple}[1]{|||#1|||}

% other macros
\newtheorem{assumption}{Assumption}
\newtheorem{theorem}{Theorem}
\newtheorem{corollary}{Corollary}
\newtheorem{lemma}{Lemma}
\newtheorem{algo}{Algorithm}
\newenvironment{remark}{\noindent{\bf Remark:}\it}{}

%>>>>>>>>>>>>>>>>>>>>>>>>>>>>>>>>>>>>>>>>>>>>>>>>>>>>>>>>>>>>>>>>>>>>>>>>>>>>>>>>>>>>>>>>>>>>>>>>>>>>>>>>>>>>>>>>>>>>>>>>>>>>>>>>>

\section{Introduction}

% general motivations, applications

Several phenomena in nature involve the interaction of a free fluid with a porous and deformable material. In geomechanics, the behavior of aquifers and groundwater flow is significantly influenced by the interaction of deformation and flow at the pore scale. In the same context, the efficient extraction of oil and gas from the subsurface relies on the ability to predict and control these phenomena. Another example from biology is the perfusion of living tissues, where the fluid carried by the main vessels is distributed by filtration to the surrounding material. All these phenomena are characterized by some common traits. In particular, they share a multiphysics nature, because they are governed by the interaction of fluid and solid mechanics.

In the general framework of continuum mechanics, several models are available to describe these phenomena. To model the free fluid, we consider for simplicity the Stokes equations, under the assumptions of incompressible and Newtonian rheology. A well accepted model for characterizing the behavior of a poroelastic material is provided by the Biot equations.  The Biot system consists of the governing equations for the deformation of an elastic skeleton, which is completely saturated with fluid. The average velocity of the fluid in the pores is modelled using the Darcy equation, complemented with an additional term that depends on the volumetric deformation of the porous matrix. Indeed, this term accounts for the poroelastic coupling. In this work we focus on the coupling of the Stokes and Biot models, for phenomena where time and space dependence of the unknowns play a significant role.

% analogies with FSI and Stokes/Darcy coupling

The numerical discretization of the problem at hand features several difficulties. With fluid-structure interaction (FSI), it shares the difficulties of combining the Eulerian description of the moving fluid domain with the typical Lagrangian parametrization of the structure. The two reference frames can be merged using the arbitrary Lagrangiann-Eulerian approach \cite{Donea1982689,Kennedy1982129}. Furthermore, loosely coupled schemes for FSI may turn out to be unconditionally unstable, under a particular range of the physical parameters of the model \cite{causin2005added,Förster20071278}. This is the so called \emph{added-mass effect}. On the other hand, accounting for a poroelastic material model requires addressing issues typical of partitioned methods for flows \cite{MR2914792,MR2970744,MR2576529} arising in particular for the coupling of Stokes-Darcy models, which represents the paradigm for studying the interaction of free flows with subsurface filtration \cite{Cao20104239,Cao20101,MR3021697,MR2899252,DMQ,Layton20022195,MR3033010,Rivière20051959}.

% review of the literature on Stokes-Biot

These observations suggest that the solution of the time dependent Biot-Stokes coupled equations is challenging from both the theoretical and computational standpoints. Concerning the analysis, the coupled problem and in particular the formulation of appropriate interface conditions has been studied in \cite{bukac2013analysis,MR2974169,MR2149168}.
From the computational viewpoint, the multiple facets of the problem have been addressed in several studies. For example, the coupling of subsurface flow and geomechanics have been recently addressed in \cite{MR2795164,MR3050000,MR2842716}. Depending on the field of application, different formulations are available to couple a free flow with a saturated poroelastic material. In the context of geosciences, this coupled problem is used to model the interaction of the material with fractures, as in \cite{ganis2,ganis1,ganis2013bem,MR2776916,MR2142590,MR2974169}. In the context of biomedical applications, FSI problems involving poroelastic materials have not been widely studied. Our work in this direction has been significantly inspired by the seminal study reported in \cite{badia2009coupling,quaini2009algorithms}.

% objective of this work

The objective of this work is to develop and analyze a loosely coupled numerical solver for the coupled Biot-Stokes system. More precisely, we design a time advancing scheme which allows us to independently solve the governing equations of the system at each time step. Resorting to time splitting approaches mitigates the difficulty to identify appropriate solvers for the coupled system and reduces the need of large memory storage. The main drawback of loosely coupled splitting schemes is possible lack of stability and accuracy. To overcome these natural limitations, we adopt a non-standard approach for the approximation of the coupling conditions, which is inspired by Nitsche's method for the enforcement of boundary conditions, and it consists of adding appropriate interface operators to the variational formulation of the problem. Using time-lagging, the variational coupled problem can be split into three independent subproblems involving the elasticity equation, Darcy equation for flow in porous media and the Stokes problem, respectively. The stability analysis of the resulting scheme shows how to design appropriate stabilization terms that guarantee the stability of the time advancing algorithm. The Nitsche's coupling approach allows for treating the mixed form of Darcy flow and thus provides accurate approximation to the filtration velocity. This is an alternative to the Lie-splitting scheme developed in \cite{bukac2013explicit}, which is suitable for the pressure formulation of Darcy flow.

% novelty w.r.t. similar approaches

A Nitsche's coupling has been proposed in \cite{MR2498525,NME:NME4607} for the interaction of a fluid with an elastic structure. Here, we extend those ideas to the case where a porous media flow is coupled to the fluid and the structure. Not only the coupling conditions between the fluid and the porous matrix subregions become significantly more complicated, but also the stability and error analysis involve the characteristic difficulties of both the FSI and Stokes-Darcy decomposition methods. The main drawback of the scheme originally proposed in \cite{MR2498525} is accuracy. In \cite{NME:NME4607} the same authors address several modifications of the scheme, mainly based on iterative corrections, to restore the optimal approximation properties that are lost when the equations of the system are split. The splitting scheme proposed in this work also affects the approximation properties of the underlying discretization methods. However, we mitigate this drawback with a different approach than in \cite{NME:NME4607}. Here, we develop a numerical solver where the loosely coupled scheme is acting as a preconditioner for the monolithic problem formulation. In this way, we combine the computational efficiency of the loosely coupled scheme with the accuracy of the discrete coupled problem. Our numerical experiments show that the number of iterations needed to solve the preconditioned problem is almost independent of the spatial and temporal discretization parameters.

% organization of the manuscript

The manuscript is organized as follows. In Section \ref{sec:problem} we present the governing equations of the prototype problem at hand, complemented by initial, boundary and interface conditions. The numerical discretization scheme, and in particular the approximation of the interface conditions is discussed in Section \ref{sec:nitsche}. A thorough analysis of the stability properties of the fully coupled discrete scheme is also presented. Section \ref{sec:uncoupled} is devoted to the development and analysis of the loosely coupled scheme. At the end of this section, we exploit the stability analysis of the method to provide insight on the properties of the loosely coupled scheme as a preconditioner for the monolithic scheme. Extensive numerical experiments are discussed in Section \ref{sec:numerics}. The corresponding results support and complement the available theory.

%>>>>>>>>>>>>>>>>>>>>>>>>>>>>>>>>>>>>>>>>>>>>>>>>>>>>>>>>>>>>>>>>>>>>>>>>>>>>>>>>>>>>>>>>>>>>>>>>>>>>>>>>>>>>>>>>>>>>>>>>>>>>>>>>>

\section{Description of the problem}\label{sec:problem}

We consider the flow of an incompressible, viscous fluid in a channel bounded by a poroelastic medium. In particular, we are interested in simulating flow through the deformable channel with a two-way coupling between the fluid and the poroleastic structure. 
The model we consider is similar to the one studied in \cite{bukac2013explicit}.
We model the flow using  the Stokes equations for a viscous, incompressible, Newtonian fluid at low speed:
\begin{align}\label{NS1}
& \rho_f \frac{\partial \boldsymbol{v}}{\partial t} = \nabla \cdot \boldsymbol\sigma_f + \boldsymbol f & \textrm{in}\; \Omega_f(t) \times (0,T), \\
 \label{NS2}
& \nabla \cdot \boldsymbol{v} = g & \textrm{in}\; \Omega_f(t) \times (0,T),
\end{align}
where $\boldsymbol{v}$ is the fluid velocity, $\boldsymbol\sigma_f= -p_f \boldsymbol{I} + 2 \mu_f \boldsymbol{D}(\boldsymbol{v})$  is the fluid stress tensor, $p_f$ is the fluid pressure, $\rho_f$ is the fluid density,  $\mu_f$ is the fluid viscosity and  $\boldsymbol{D}(\boldsymbol{v}) = (\nabla \boldsymbol{v}+(\nabla \boldsymbol{v})^{T})/2$ is the rate-of-strain tensor.

\begin{figure}
\begin{center}
\begin{picture}(0,0)%
\includegraphics{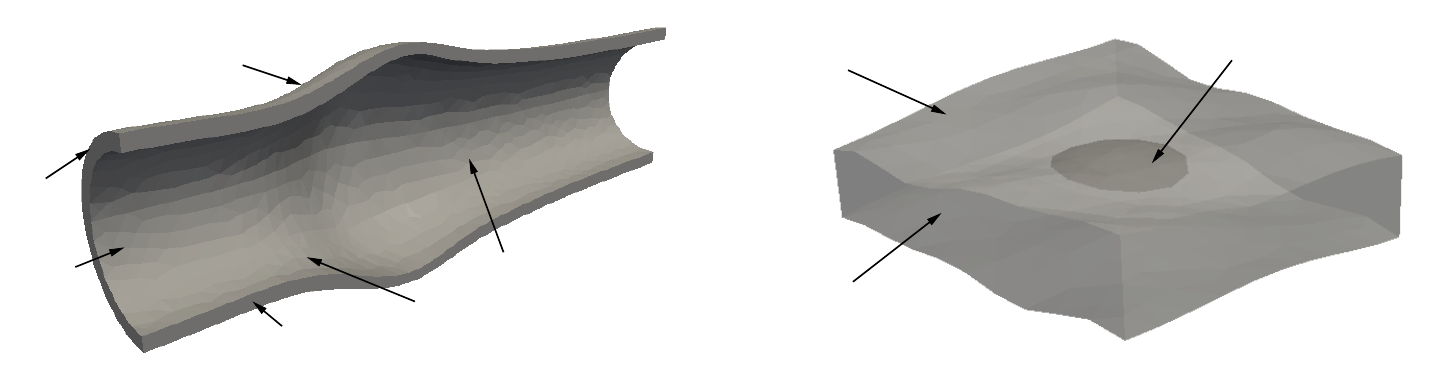}%
\end{picture}%
\setlength{\unitlength}{2072sp}%
\begingroup\makeatletter\ifx\SetFigFont\undefined%
\gdef\SetFigFont#1#2#3#4#5{%
  \reset@font\fontsize{#1}{#2pt}%
  \fontfamily{#3}\fontseries{#4}\fontshape{#5}%
  \selectfont}%
\fi\endgroup%
\begin{picture}(13227,3450)(571,-3286)
\put(3196,-3031){\makebox(0,0)[lb]{\smash{{\SetFigFont{10}{12.0}{\familydefault}{\mddefault}{\updefault}{\color[rgb]{0,0,0}$\Omega_p$}%
}}}}
\put(4366,-2851){\makebox(0,0)[lb]{\smash{{\SetFigFont{10}{12.0}{\familydefault}{\mddefault}{\updefault}{\color[rgb]{0,0,0}$\Gamma$}%
}}}}
\put(5221,-2401){\makebox(0,0)[lb]{\smash{{\SetFigFont{10}{12.0}{\familydefault}{\mddefault}{\updefault}{\color[rgb]{0,0,0}$\Omega_f$}%
}}}}
\put(901,-2581){\makebox(0,0)[lb]{\smash{{\SetFigFont{10}{12.0}{\familydefault}{\mddefault}{\updefault}{\color[rgb]{0,0,0}$\Gamma_f^{in}$}%
}}}}
\put(2521,-241){\makebox(0,0)[lb]{\smash{{\SetFigFont{10}{12.0}{\familydefault}{\mddefault}{\updefault}{\color[rgb]{0,0,0}$\Gamma_p^{ext}$}%
}}}}
\put(6346,-826){\makebox(0,0)[lb]{\smash{{\SetFigFont{10}{12.0}{\familydefault}{\mddefault}{\updefault}{\color[rgb]{0,0,0}$\Gamma_f^{out}$}%
}}}}
\put(6706,-1321){\makebox(0,0)[lb]{\smash{{\SetFigFont{10}{12.0}{\familydefault}{\mddefault}{\updefault}{\color[rgb]{0,0,0}$\Gamma_p^{out}$}%
}}}}
\put(8236,-241){\makebox(0,0)[lb]{\smash{{\SetFigFont{10}{12.0}{\familydefault}{\mddefault}{\updefault}{\color[rgb]{0,0,0}$\Gamma_p^{ext}$}%
}}}}
\put(8236,-2716){\makebox(0,0)[lb]{\smash{{\SetFigFont{10}{12.0}{\familydefault}{\mddefault}{\updefault}{\color[rgb]{0,0,0}$\Omega_p$}%
}}}}
\put(11791,-196){\makebox(0,0)[lb]{\smash{{\SetFigFont{10}{12.0}{\familydefault}{\mddefault}{\updefault}{\color[rgb]{0,0,0}$\Omega_f,\,\Gamma$}%
}}}}
\put(586,-1816){\makebox(0,0)[lb]{\smash{{\SetFigFont{10}{12.0}{\familydefault}{\mddefault}{\updefault}{\color[rgb]{0,0,0}$\Gamma_p^{in}$}%
}}}}
\end{picture}%

\end{center}
\caption{Description of the domain and boundary decomposition for the two problem configurations addressed in this work, an arterial segment (left) and a fractured reservoir (right). In the latter case, the fracture has a discoidal shape and it is completely embedded into the reservoir.}\label{fig:domains}
\end{figure}
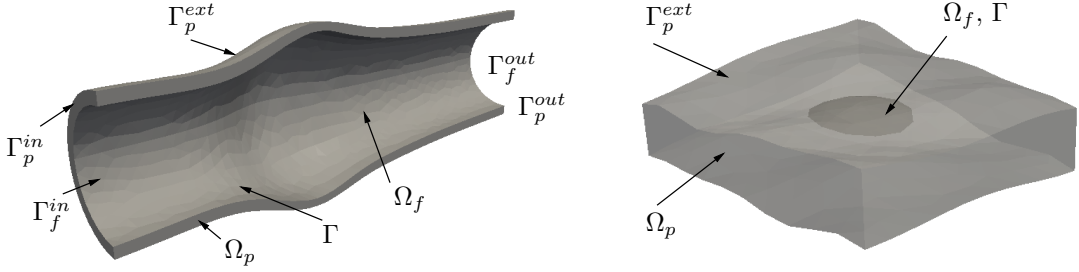

We consider two configurations: ($i$) the channel extends to the external boundary and ($ii$) the channel is surrounded by the poroelastic media, see Figure \ref{fig:domains}. Configuration ($i$) is suitable for FSI in arteries, while case ($ii$) applies to fractured reservoirs. We present the method for configuration ($i$), which has more complex boundary conditions. In this case, denote the inlet and outlet fluid boundaries by $\Gamma^{in}_f$ and $\Gamma_f^{out}$, respectively. At the inlet and outlet boundary we prescribe the following conditions:
\begin{align}
& \boldsymbol {\sigma}_f \boldsymbol {n}_{f} = - p_{in}(t) \boldsymbol{n}_{f} 
\quad \text{or} \quad \boldsymbol v = \boldsymbol v_{in}(t)
& \textrm{on} \; \Gamma_f^{in} \times (0,T), \label{inlet} \\
& \boldsymbol {\sigma}_f \boldsymbol {n}_{f} = 0  & \textrm{on} \; \Gamma_f^{out} \times (0,T),\label{outlet}
\end{align}
where $\boldsymbol{n}_{f}$ is the outward normal unit vector to the fluid boundaries and $p_{in}(t)$ is the pressure increment with respect to the ambient pressure surrounding the channel.

The fluid domain is bounded by a deformable porous matrix consisting of a skeleton and connecting pores filled with fluid, whose dynamics is described by the Biot model, which in the Eulerian formulation reads as follows:
\begin{align}
& \rho_{p} \frac{D^2 \boldsymbol U}{D t^2} - \nabla \cdot \boldsymbol \sigma_p = \boldsymbol h & \textrm{in} \; \Omega_p(t)\times  (0,T),  \label{B1}\\
& {\kappa}^{-1} \boldsymbol q = -\nabla p_p & \textrm{in} \; \Omega_p(t) \times (0,T),  \label{B2} \\
& \frac{D}{D t}(s_0 p_p + \alpha \nabla \cdot \boldsymbol U) + \nabla \cdot \boldsymbol q = s & \textrm{in} \; \Omega_p(t) \times  (0,T),\label{B3}
\end{align}
where $\frac{D}{Dt}$ denotes the material derivative. The stress tensor of the poroelastic medium is given by
$
\boldsymbol \sigma_p = \boldsymbol \sigma^E - \alpha p_p \boldsymbol I,
$
where $\boldsymbol \sigma^E$ denotes the elasticity stress tensor. With the assumption that the displacement $\boldsymbol U$ of the skeleton is connected to stress tensor $\boldsymbol \sigma^E$ via the Saint Venant-Kirchhoff elastic model, we have
$
\boldsymbol \sigma^E (\boldsymbol U) = 2 \mu_p  \boldsymbol D (\boldsymbol U) + \lambda_p \textrm{tr}(\boldsymbol D(\boldsymbol U)) \boldsymbol I,
$
where $\lambda_p$ and $\mu_p$ denote the Lam\'e coefficients for the skeleton, and, with the hypothesis of ``small'' deformations, $\boldsymbol D(\boldsymbol U) = (\nabla \boldsymbol U+(\nabla \boldsymbol U)^{T})/2.$ 
System~\eqref{B1}-\eqref{B3} consists of the momentum equation for the balance of total forces~\eqref{B1}, Darcy's law~\eqref{B2}, and the storage equation~\eqref{B3} for the fluid mass conservation in the pores of the matrix, where the flux $\boldsymbol q$ is the relative velocity of the fluid within the porous structure and $p_p$ is the fluid pressure. The density of saturated porous medium is denoted by $\rho_{p}$, and the hydraulic conductivity is denoted by $\kappa$. To account for anisotropic transport properties, $\kappa$ is in general a symmetric positive definite tensor. For simplicity of notation, but without loss of generality with respect to the derivation of the numerical scheme, in what follows we consider it as a scalar quantity. The coefficient $s_0 \in (0,1)$ is the storage coefficient, and the Biot-Willis constant $\alpha$ is the pressure-storage coupling coefficient. 

For configuration ($i$), we assume that the poroelastic structure is fixed at the inlet and outlet boundaries:
\begin{equation}\label{homostructure1}
\boldsymbol U = \boldsymbol 0 \quad \text{on} \ \Gamma_p^{in} \cup \Gamma_p^{out} \times (0,T),
\end{equation}
that the external structure boundary $\Gamma_p^{ext}$ is exposed to external ambient pressure 
\begin{eqnarray}
 \boldsymbol {n}_{p} \cdot \boldsymbol \sigma^E \boldsymbol {n}_{p} &=&  0 \quad \textrm{on} \; \Gamma_p^{ext} \times (0,T),
\end{eqnarray}
where  $\boldsymbol n_{p}$ is the outward unit normal vector on $\partial\Omega_p$,
and that the tangential displacement of the exterior boundary is zero:
\begin{equation}
\boldsymbol U \cdot \ttau_p =  0 \quad \textrm{on} \; \Gamma_p^{ext} \times (0,T),
\end{equation}
where $U \cdot \ttau_p$ denotes the tangential component of the vector $\boldsymbol U$.
On the fluid pressure in the porous medium, we impose following boundary conditions:
\begin{equation}
p_p = 0 \quad \textrm{on} \; \Gamma_p^{ext},
\quad
\boldsymbol q \cdot \boldsymbol n_p = 0 \quad \textrm{on} \; \Gamma_p^{in} \cup \Gamma_p^{out}  \times (0,T).
\end{equation}
For configuration ($ii$) we assume $\boldsymbol U = 0$ and $\boldsymbol q \cdot \boldsymbol n_p = 0$ on $\partial \domp$.

At the initial time, the fluid and the poroelastic structure are assumed to be at rest, with zero displacement from the reference configuration
\begin{equation}\label{initial}
\boldsymbol{v}=0, \quad \boldsymbol U = 0, \quad \frac{D \boldsymbol U}{D t}=0, \quad \boldsymbol q =0.
\end{equation}

Finally, the fluid and poroelastic structure are coupled via the following interface conditions, where we denote by $\boldsymbol n$ the outward normal to the fluid domain and by $\boldsymbol t$ the tangential unit vector on the interface $\Gamma(t)$. We assume that $\boldsymbol n,\,\boldsymbol t$ coincide with the unit vectors relative to the fluid domain $\domf$.
For mass conservation, the continuity of normal flux implies that
\begin{equation}\label{CNF}
\Big( \boldsymbol{v} - \aalpha \frac{D \boldsymbol U}{D t} \Big)  \cdot \boldsymbol n  =  \boldsymbol q \cdot \boldsymbol n  
\quad \textrm{on} \; \Gamma(t).
\end{equation} 
There are different options to formulate a condition relative for the tangential velocity field at the interface. A no-slip interface condition is appropriate for those problems where fluid flow in the tangential direction is not allowed,
\begin{equation}\label{CBJS}
\boldsymbol v \cdot \boldsymbol t = \frac{D\boldsymbol U}{D t} \cdot \boldsymbol t \quad \textrm{on} \; \Gamma(t).
\end{equation}
In alternative, a Beavers-Joseph-Saffman type condition may be prescribed \cite{Layton20022195,MR2149168},
\begin{equation}\label{CBJS_bis}
\beta \Big( \boldsymbol v - \frac{D\boldsymbol U}{D t} \Big) \cdot \boldsymbol t 
= - \boldsymbol t \cdot \boldsymbol{\sigma}_f \boldsymbol{ n}
\quad \textrm{on} \; \Gamma(t),
\end{equation}
where the parameter $\beta$ quantifies the resistance that the porous matrix opposes to fluid flow in the tangential direction. Concerning the exchange of stresses, the balance of normal components of the stress in the fluid phase gives:
\begin{equation}\label{CBNSF}
 \boldsymbol n \cdot \boldsymbol \sigma_f \boldsymbol n = -p_p  \quad \textrm{on} \; \Gamma(t).
\end{equation}
The conservation of momentum describes balance of contact forces. 
Precisely, it says that the sum of contact forces at the fluid-porous medium interface is equal to zero:
\begin{eqnarray}
\label{CBFN} 
\aalpha \boldsymbol n \cdot \boldsymbol{\sigma}_f \boldsymbol n - \boldsymbol n \cdot \boldsymbol {\sigma}_p \boldsymbol n = 0 
\quad \textrm{on} \; \Gamma(t),
\\
\label{CBFT}
\boldsymbol t \cdot \boldsymbol{\sigma}_f \boldsymbol{ n} - \boldsymbol t \cdot \boldsymbol {\sigma}_p \boldsymbol n = 0 
\quad \textrm{on} \; \Gamma(t). 
\end{eqnarray}

%>>>>>>>>>>>>>>>>>>>>>>>>>>>>>>>>>>>>>>>>>>>>>>>>>>>>>>>>>>>>>>>>>>>>>>>>>>>>>>>>>>>>>>>>>>>>>>>>>>>>>>>>>>>>>>>>>>>>>>>>>>>>>>>>>

\section{Finite element approximation of the coupled problem using Nitsche's method for the interface conditions}\label{sec:nitsche}

For the sake of simplicity, we present and analyze the FSI algorithms under the assumption of fixed domains $\Omega_f(t)$ and $\Omega_p(t)$, namely $$\Omega_f(t)=\Omega_f, \quad \Omega_p(t)=\Omega_p, \quad \forall t \in (0,T).$$ Although simplified, this problem still retains the main difficulties associated with the added-mass effect and with the fluid-porous media coupling. 
Nitsche's method is a technique to weakly enforce coupling conditions at the discrete level \cite{hansbo2005nitsche}. For the spatial discretization, we exploit the finite element method. We denote with $\femvf,\fempf$ the finite element spaces for the velocity and pressure approximation on the fluid domain $\domf$, with $\femvp,\fempp$ the spaces for velocity and pressure approximation on the porous matrix $\domp$ and with $\fems,\dotfems$ the approximation spaces for the structure displacement and velocity, respectively.   We assume that all the finite element approximation spaces comply with the prescribed Dirichlet conditions on external boundaries $\partial\domf,\,\partial\domp$, without entering into the details of all possible variations for cases ($i$) and ($ii$).
The bilinear forms relative to the structure, namely equation \eqref{B1}, are:
\begin{align*}
\bilinasp(\dispp,\bsp) :=& 
2 \mu_p \int_{\Omega_p} \boldsymbol D(\dispp) : \boldsymbol D(\bsp) d \boldsymbol x 
+ \lambda_p \int_{\Omega_p} (\nabla \cdot \dispp)(\nabla \cdot \bsp)  d \boldsymbol x,
\\
\bilinbsp(\pp,\bsp) :=& \alpha \int_{\domp} \pp \nabla \cdot \bsp d \boldsymbol x.
\end{align*}

\noindent
For the free flow, equations \eqref{NS1}-\eqref{NS2}, and the filtration through the porous matrix, equations \eqref{B2}-\eqref{B3}, the bilinear forms are:
\begin{equation*}
\begin{array}{ll}
\displaystyle
\bilinaff(\velf,\bff) := 2 \mu_f \int_{\domf} \boldsymbol D(\velf) : \boldsymbol D(\bff) d \boldsymbol x,
& 
\displaystyle
\bilinafp(\velp,\bfp) := \int_{\domp} \kappa^{-1} \velp \cdot \bfp d \boldsymbol x,
\\[10pt]
\displaystyle
\bilinbff(\pf,\bff) := \int_{\domf} \pf \nabla \cdot \bff d \boldsymbol x,
&
\displaystyle
\bilinbfp(\pp,\bfp) := \int_{\domp} \pp \nabla \cdot \bfp d \boldsymbol x .
\end{array}
\end{equation*}
After integrating by parts the governing equations, in order to distribute over test functions the second spatial derivatives of velocities and displacements as well as the first derivatives of the pressure, resorting to the \emph{dual-mixed} weak formulation of Darcy's problem, the following interface terms appear in the variational equations,
\begin{equation*}
I_{\interf} =  \int_{\interf} (\strf \nn \cdot \bff - \strp \nn \cdot \bsp  +  \pp \bfp \cdot \nn  ).
\end{equation*}
Starting from the expression of $I_{\interf}$, Nitsche's method allows us to weakly enforce the interface conditions \eqref{CNF}-\eqref{CBFT}. More precisely, we separate $I_{\interf}$ into the normal and tangential components with respect to $\Gamma$ and we use \eqref{CBNSF}, \eqref{CBFN}, \eqref{CBFT} to substitute the components of $\strf$ into $\strp$ and $\pp$. As a result $I_{\interf}$ can be rewritten as,
\begin{equation*}
I_{\interf} = \int_{\interf} \nn \cdot \strf(\velf,\pf) \nn\, \big( \bff - \bfp - \aalpha \bsp \big)\cdot \nn
+\int_{\interf} \ttau \cdot \strf(\velf,\pf) \nn (\bff - \bsp ) \cdot \ttau.
\end{equation*}
Since the expression of the interface terms involves the stresses only on the fluid side, this formulation can be classified as a \emph{one-sided} variant of Nitsche's method for interface conditions. We refer to \cite{hansbo2005nitsche} for an overview of different formulations.
The enforcement of the kinematic conditions (\eqref{CNF} and \eqref{CBJS}) using Nitsche's method is based on adding to the variational formulation of the problem appropriate penalty terms. This results in the transformed integral,
\begin{align*}
&- I_{\interf}^*(\velf,\velp,\pf,\pp,\dispp; \bff,\bfp,\bpf,\bpp,\bsp) =
\\
& - \int_{\interf} \nn \cdot \strf(\velf,\pf) \nn\, \big( \bff - \bfp - \aalpha \bsp \big)\cdot \nn
  - \int_{\interf} \ttau \cdot \strf(\velf,\pf) \nn (\bff - \bsp ) \cdot \ttau
\\
& + \int_{\interf} \gamma_f \mu_f h^{-1} \big[ \big( \velf - \velp - \aalpha \partial_t \dispp \big)\cdot \nn\, \big( \bff - \bfp - \aalpha \bsp \big)\cdot \nn
   + \big(\velf - \partial_t \dispp \big) \cdot \ttau\, \big(\bff -  \bsp \big) \cdot \ttau \big],
\end{align*}
where $\gamma_f>0$ denotes a penalty parameter that will be suitably defined later on.
Furthermore, in order to account for the symmetric, incomplete or skew-symmetric variants of Nitsche's method, see \cite{hansbo2005nitsche}, we introduce the following additional terms
\begin{align*}
&- S_{\interf}^{*,\sym}(\velf,\velp,\pf,\pp,\dispp; \bff,\bfp,\bpf,\bpp,\bsp) =
\\
& - \int_{\interf} \nn \cdot \strf(\sym\bff,-\bpf) \nn\, \big( \velf - \velp - \aalpha \partial_t \dispp \big)\cdot \nn
  - \int_{\interf} \ttau \cdot \strf(\sym\bff,-\bpf) \nn (\velf- \partial_t\dispp ) \cdot \ttau,
\end{align*}
which anyway do not violate the consistency of the original scheme because they vanish if the kinematic constraints are satisfied. The flag $\sym \in \{1,0,-1\}$ determines if we adopt a symmetric, incomplete or skew symmetric formulation, respectively. We can also easily accommodate the weak enforcement of the Beavers-Joseph-Saffman condition, namely \eqref{CBJS_bis}. In this case, the interface operators become,
\begin{align*}
I_{\interf}^++S_{\interf}^{+,\sym} = & - \int_{\interf} \nn \cdot \strf(\velf,\pf) \nn\, \big( \bff - \bfp - \aalpha \bsp \big)\cdot \nn
\\
& - \int_{\interf} \nn \cdot \strf(\sym\bff,-\bpf) \nn\, \big( \velf - \velp - \aalpha \partial_t \dispp \big)\cdot \nn
\\
& + \int_{\interf} \gamma_f \mu_f h^{-1} \big( \velf - \velp - \aalpha \partial_t \dispp \big)\cdot \nn\, \big( \bff - \bfp - \aalpha \bsp \big)\cdot \nn
\\
& + \int_{\interf} \beta \big(\velf - \partial_t \dispp \big) \cdot \ttau\, \big(\bff -  \bsp \big) \cdot \ttau.
\end{align*}
By comparing $I_{\interf}^+ + S_{\interf}^{+,\sym}$ and $I_{\interf}^* + S_{\interf}^{*,\sym}$ we observe that the operators corresponding to \eqref{CBJS_bis} can be seen as a particular form of the more general case that is obtained when no-slip conditions \eqref{CBJS} are enforced weakly, namely $I_{\interf}^* + S_{\interf}^{*,\sym}$. For this reason, we perform the analysis of the numerical scheme in the latter form.

For any $t \in (0,T)$, the coupled fluid/solid problem consists of finding
$\velf,\pf,\velp,\pp \in \femvf \times \fempf \times \femvp \times \fempp$ 
and $\dispp,\dotdispp\in \fems \times \dotfems$ such that for any 
$\bff,\bpf,\bfp,\bpp \in \femvf \times \fempf \times \femvp \times \fempp$ 
and  $\bsp,\dotbsp \in \fems \times \dotfems$ we have,
\begin{multline}\label{eq:SDP}
\rho_p \int_{\domp} \partial_t \dotdispp \cdot \bsp  d\boldsymbol x
+ \rho_p \int_{\domp} \big(\dotdispp  - \partial_t \dispp \big) \cdot \dotbsp d\boldsymbol x
+ \rho_f \int_{\domf} \partial_t \velf \cdot \bff d\boldsymbol x
+ s_0 \int_{\domp} \partial_t \pp \bpp d\boldsymbol x
\\
+ \bilinasp(\dispp,\bsp) - \bilinbsp(\pp,\bsp) + \bilinbsp(\bpp,\partial_t\dispp) 
+ \bilinafp(\velp,\bfp) - \bilinbfp(\pp,\bfp) + \bilinbfp(\bpp,\velp) 
\\
+ \bilinaff(\velf,\bff) - \bilinbff(\pf,\bff) + \bilinbff(\bpf,\velf) 
 - (I_{\interf}^*+S_{\interf}^{*,\sym})(\velf,\velp,\pf,\pp,\dispp; \bff,\bfp,\bpf,\bpp,\bsp) 
\\
= \fff(t;\bff,\bfp,\bpf,\bpp,\bsp).
\end{multline}
Problem \eqref{eq:SDP} is usually called the semi-discrete problem (SDP). For simplicity of notation, the bilinear form corresponding to the left hand side of \eqref{eq:SDP} will be denoted as $\aas(\cdot,\cdot)$. Equation \eqref{eq:SDP} must be complemented by suitable initial conditions and $\fff(\cdot)$ accounts for boundary conditions and forcing terms, which vary from case ($i$) to case ($ii$). For example, for case ($i$) we set $\boldsymbol f=\boldsymbol h=\boldsymbol 0$, $g=s=0$ and $p_{in} \neq 0$ or $\boldsymbol v_{in} \neq \boldsymbol 0$ on $\Gamma_f^{in}$. For case ($ii$) we have $\boldsymbol f=\boldsymbol h=\boldsymbol 0$, $s=0$, but we consider a non vanishing flow source $g \neq \boldsymbol 0$. The corresponding forcing terms are,
\begin{equation*}
(i) \quad \fff(t;\bff) = - \int_{\Gamma_f^{in}} p_{in}(t)\, \bff \cdot \nn_f,
\quad
(ii) \quad \fff(t;\bff) = \int_{\domf} g(\boldsymbol x,t)\, \bpf.
\end{equation*}
To avoid addressing the ramifications relative to these particular cases in the forthcoming analysis, we denote by $\fff(t;\bff)$ a generic forcing term acting on the flow equations and we assume there exists a norm $\|\fff(t_n)\|$ such that
\begin{equation}\label{eq:fff}
\fff(t;\bff) \leq \|\fff(t_n)\| \|\boldsymbol D(\bff)\|_{\domf},
\end{equation}
where from now on $\|\cdot\|_D$ denotes the $L^2$-norm on the set $D$.

We now address the time discretization. We denote with $t_n$ the current time step and with $\dt$ the first order (backward) discrete time derivative $\dt u^n := \tau^{-1} (u^n-u^{n-1})$. The fully discrete coupled fluid / solid problem consists of finding, at each time step $t_n$,
$\velf^n,\pf^n,\velp^n,\pp^n \in \femvf \times \fempf \times \femvp \times \fempp$ 
and $\dispp^n,\dotdispp^n\in \fems \times \dotfems$ such that for any 
$\bff,\bpf,\bfp,\bpp \in \femvf \times \fempf \times \femvp \times \fempp$ 
and  $\bsp,\dotbsp \in \fems \times \dotfems$ we have,
\begin{multline}\label{eq:DIP}
\rho_p \int_{\domp} \dt \dotdispp^n \cdot \bsp  d\boldsymbol x
+ \rho_p \int_{\domp} \big(\dotdispp^n  - \dt \dispp^n \big) \cdot \dotbsp d\boldsymbol x
+ \rho_f \int_{\domf} \dt \velf^n \cdot \bff d\boldsymbol x
+ s_0 \int_{\domp} \dt \pp^n \bpp d\boldsymbol x
\\
+ \bilinasp(\dispp^n,\bsp) - \bilinbsp(\pp^n,\bsp) + \bilinbsp(\bpp,\dt\dispp^n) 
+ \bilinafp(\velp^n,\bfp) - \bilinbfp(\pp^n,\bfp) + \bilinbfp(\bpp,\velp^n) 
\\
+ \bilinaff(\velf^n,\bff) - \bilinbff(\pf^n,\bff) + \bilinbff(\bpf,\velf^n) 
 - (I_{\interf}^*+S_{\interf}^{*\sym})(\velf^n,\velp^n,\pf^n,\pp^n,\dispp^n; \bff,\bfp,\bpf,\bpp,\bsp) 
\\
= \fff(t^n;\bff).
\end{multline}
We denote problem \eqref{eq:DIP} as the (fully) discrete problem with implicit coupling between the fluid and the structure sub-problems (DIP). The bilinear form on the left hand side is represented with the symbol $\aasi(\cdot,\cdot)$. We note that for the time discretization of \eqref{eq:DIP} we have adopted the Backward Euler (BE) method for both the flow and the structure problem. 

\begin{remark}
The structure problem may be also discretized using the Newmark scheme, which enjoys better energy conservation properties than BE \cite{MR1299729}. Switching to the Newmark scheme would not affect the discretization of the interface conditions. In other words, their approximation using Nitsche's method is completely independent on the time discretization scheme. However, later on we will point out that the numerical dissipation artificially introduced when using the BE method improves the stability of the scheme when loosely coupled splitting strategies are considered.
\end{remark}

%--------------------------------------------------------------------------------------------------------
The next step consists of analyzing the stability of the fully coupled problem. Before proceeding, we define the energy of the system,
$\energyf^n$ and $\energys^n$ for the fluid and the poroelastic structure respectively, as follows:
\begin{align*}
\energyf^n &:= 
\frac{1}{2} \rho_f \|\velf^n\|^2_{\domf},
\\
\energys^n &:= 
\frac12 \Big( \rho_p\| \dotdispp^n \|^2_{\domp}   
+ 2 \mu_p \|{\boldsymbol D}(\dispp^n)\|^2_{\domp}
+ \lambda_p \|\nabla \cdot \dispp^n\|^2_{\domp} 
+ s_0\|\pp^n\|^2_{L^2(\domp)} \Big).
\end{align*}
%--------------------------------------------------------------------------------------------------------
The following algebraic identity will be systematically used in the forthcoming derivations,
\begin{equation}\label{eq:dtt}
\int_{\Omega} u^n \dt u^n = \frac12 \dt \| u^n \|_{\Omega}^2 + \frac12 \tau \| \dt u^n \|_{\Omega}^2.
\end{equation}
We will also use the following inverse inequality, which holds true for a family of shape-regular and quasi-uniform meshes \cite{MR2050138}, where $C_{TI}$ is a positive constant independent of the mesh characteristic size $h$,
\begin{equation}\label{eq:trace_inverse}
h \|\DD(\uu_h) \nn\|_{\interf}^2 \leq C_{TI} \|\DD(\uu_h)\|_{\domf}^2.
\end{equation}
The following result shows that the time advancing scheme used for the coupled problem is stable under conditions on the penalty parameters $\gamma_f$.

\begin{theorem}\label{th:energy}
For any $\epsfa,\epsfab$ that satisfy
$$\big(1 - \frac{(\sym+1)}{2}\epsfa C_{TI} -\frac{\epsfab}{2}  \big)>0$$
where $\sym \in\{-1,0,1\}$ provided that $\gamma_f > (\sym +1) (\epsfa)^{-1}$,
there exist constants $0<c<1$ and $C>1$, uniformly independent of the mesh characteristic size $h$, such that
\begin{align}
\label{eq:stability-DIP}
& \energyf^N + \energys^N 
+ c \tau \sum_{n=1}^N \Big[ 2 \mu_f \|\DD(\velf^n)\|^2_{\domf} + \kappa^{-1} \|\velp^n\|^2_{\domp} 
\\
\nonumber
& + \frac{\tau}{2} \big( \rho_f \|\dt \velf^n\|_{\domf}^2 
+ 2 \mu_p \|\dt \DD(\dispp^n)\|^2_{\domp} 
+ s_0 \|\dt \pp\|_{\domp}^2
+ \lambda_p \| \dt \nabla \cdot \dispp \|_{\domp}^2 \big)
\\
\nonumber
& + \mu_f h^{-1} \big( \| \big( \velf^n - \velp^n - \aalpha \dt \dispp^n \big)\cdot \nn \|_{\interf}^2
+ \| \big(\velf^n - \dt \dispp^n \big) \cdot \ttau \|_{\interf}^2 \big)
\Big]
\\
\nonumber
& \leq \energyf^0 + \energys^0 
+ \tau \sum_{n=1}^N \frac{C}{\mu_f} \|\fff(t_n)\|^2.
\end{align}
More precisely, we have
\begin{align*}
c &< \min \{ \big(1 -\frac{(\sym+1)}{2}\epsfa C_{TI} -\frac{\epsfab}{2} \big), \big((\gamma_f - (\sym+1)(\epsilon_f^\prime)^{-1}\big) \},
\\
C &> (2 \epsfab)^{-1}.
\end{align*}
\end{theorem}
%--------------------------------------------------------------------------------------------------------
{\bf Proof:}
We test \eqref{eq:DIP} with the following functions,
\begin{equation*}
\bff=\velf^n,\,\bfp=\velp^n,\,\bpf=\pf^n,\,\bpp=\pp^n,\,\bsp=\dt \dispp^n,\, \dotbsp=\dt \dotdispp^n.
\end{equation*}

%--------------------------------------------------------------------------------------------------------
For the terms related to the flow problem we have,
\begin{multline*}
\int_{\domf} \rho_f \dt \velf^n \cdot \velf^n 
+ \bilinaff(\velf^n,\velf^n) - \bilinbff(\pf^n,\velf^n) + \bilinbff(\pf^n,\velf^n) 
+ \bilinafp(\velp^n,\velp^n) - \bilinbfp(\pp^n,\velp^n) + \bilinbfp(\pp^n,\velp^n) 
\\
= \dt \energyf + \frac{\rho_f}{2} \tau \|\dt \velf^n\|^2_{\domf} + 2 \mu_f \|\DD(\velf^n)\|^2_{\domf} 
+ \kappa^{-1} \|\velp^n\|^2_{\domp} .
\end{multline*}
Furthermore, the interface terms for the coupling between the fluid and the structure can be bounded as follows,
\begin{multline*}
- \int_{\interf} \nn \cdot \strf(\velf^n,\pf^n) \nn\, \big( \velf^n - \velp^n - \aalpha \dt \dispp^n \big)\cdot \nn
- \int_{\interf} \nn \cdot \strf(\sym \velf^n,-\pf^n) \nn\, \big( \velf^n - \velp^n - \aalpha \dt \dispp^n \big)\cdot \nn
\\
= - (1 + \sym) \int_{\interf} \nn \cdot \big(2\mu_f\DD(\velf^n)\big) \nn\, \big( \velf^n - \velp^n - \aalpha \dt \dispp^n \big)\cdot \nn
\leq 2\mu_f (1 + \sym) \|\DD(\velf^n)\nn\|_\interf \|\big( \velf^n - \velp^n - \aalpha \dt \dispp^n \big)\cdot \nn\|_\interf
\\
\leq \mu_f (1 + \sym) \epsfa C_{TI} \|\DD(\velf^n)\nn\|_{\domf}^2 
+ \mu_f (1 + \sym) (\epsfa h)^{-1} \|\big( \velf^n - \velp^n - \aalpha \dt \dispp^n \big)\cdot \nn\|_\interf^2
\end{multline*}
The interface terms containing the tangential component of the fluid stress, 
namely $\ttau \cdot \strf(\velf^n,\pf^n) \nn$ can be estimated analogously.

%--------------------------------------------------------------------------------------------------------
For the structure problem we have,
\begin{multline*}
\rho_p \int_{\domp} \dt \dotdispp^n \cdot \dt \dispp^n
+ \rho_p \int_{\domp} \big(\dotdispp^n - \dt \dispp^n \big) \cdot \dt \dotdispp^n
\\
+ \bilinasp(\dispp^n,\dt\dispp^n) 
- \bilinbsp(\pp^n,\dt\dispp^n) 
+ \bilinbsp(\pp^n,\dt\dispp^n)
+ \bilincfp(\pp^n,\pp^n)
\\
\leq \dt \energys 
+ \frac{\tau}{2} \big( 2 \mu_p \|\dt{\boldsymbol D}(\dispp^n)\|^2_{\domp} 
+ \lambda_p \| \dt \nabla \cdot \dispp^n \|_{\domp}^2 
+ s_0 \|\dt \pp^n\|_{\domp}^2 \big)
\end{multline*}

%--------------------------------------------------------------------------------------------------------
Lastly, owing to \eqref{eq:fff} we have 
\begin{equation*}
\fff(t_n;\velf^n)
\leq (2\epsfab \mu_f)^{-1} \|\fff(t_n)\|^2
+ \frac{\epsfab}{2} \mu_f \|\DD(\velf^n)\|_{\domf}^2.
\end{equation*}

%--------------------------------------------------------------------------------------------------------
To conclude, we sum up with respect to the time index $n=1,\ldots,N$ and we multiply by $\tau$.
As a result of the previous inequalities, for any $\epsfa,\epsfab>0$,  
the following energy estimate holds true,
\begin{align}
\nonumber
& \energyf^N +  \energys^N 
+ \tau \sum_{n=1}^N \Big[ 2\mu_f \big(1 - \frac{(\sym+1)}{2}\epsfa C_{TI} -\frac{\epsfab}{4}  \big) \|\DD(\velf^n)\|^2_{\domf} 
+ \kappa^{-1} \|\velp^n\|^2_{\domp} 
\\
\nonumber
& + \frac{\tau}{2} \big( \rho_f \|\dt \velf^n\|_{\domf}^2
   + 2 \mu_p \|\dt\DD(\dispp^n)\|^2_{\domp}
   + s_0 \|\dt \pp^n\|_{\domp}^2 
   + \lambda_p \| \dt \nabla \cdot \dispp^n \|_{\domp}^2 \big)
\\
\nonumber
& +  \big((\gamma_f - (\sym+1)(\epsilon_f^\prime)^{-1}\big) \mu_f h^{-1} \big( \| \big( \velf^n - \velp^n - \aalpha \dt \dispp^n \big)\cdot \nn \|_{\interf}^2
+ \| \big(\velf^n - \dt \dispp^n \big) \cdot \ttau \|_{\interf}^2 \big)
\Big]
\\
\label{eq:energy-coupled}
& \leq \energyf^0 + \energys^0 
 + \tau \sum_{n=1}^N (2\epsfab \mu_f)^{-1} \|\fff(t_n)\|^2,
\end{align}
which implies \eqref{eq:stability-DIP}.
\hfill$\Box$

%>>>>>>>>>>>>>>>>>>>>>>>>>>>>>>>>>>>>>>>>>>>>>>>>>>>>>>>>>>>>>>>>>>>>>>>>>>>>>>>>>>>>>>>>>>>>>>>>>>>>>>>>>>>>>>>>>>>>>>>>>>>>>>>>>

%>>>>>>>>>>>>>>>>>>>>>>>>>>>>>>>>>>>>>>>>>>>>>>>>>>>>>>>>>>>>>>>>>>>>>>>>>>>>>>>>>>>>>>>>>>>>>>>>>>>>>>>>>>>>>>>>>>>>>>>>>>>>>>>>>
\section{A loosely coupled time advancing scheme}\label{sec:uncoupled}

When enforced by Nitsche's method, the interface conditions \eqref{CNF}-\eqref{CBFT} appear in the variational formulation in a \emph{modular} form. As a result, using time lagging, it is straightforward to design various loosely coupled algorithms to solve each equation of the problem independently from the others. The main issue is to identify those algorithms which enjoy the best stability properties. Here, we address two interesting possibilities. The first \emph{loosely coupled} or explicit time advancing scheme can be obtained by solving at each time step $t^n$ the following sub-problems.

\medskip

\noindent{\bf Algorithm (a):}
\begin{description}
\item[$(a.1)$] Solve the \emph{Biot problem}: writing the interface terms explicitly and using the unit normal and tangential vectors relative to $\domp$, given $\velf^{n-1},\pf^{n-1} \in \femvf \times \fempf$  we aim to find  $\velp^n,\pp^n,\dispp^n,\dotdispp^n \in \femvp \times \fempp \times \fems \times \dotfems$ such that
\begin{align}
\nonumber
& \rho_p \int_{\domp} \dt \dotdispp^n \cdot \bsp d \boldsymbol x
+ \rho_p \int_{\domp} \big( \dotdispp^n - \dt \dispp^n \big) \cdot \dotbsp d \boldsymbol x
+ s_0 \int_{\domp} \dt \pp^n \bpp d \boldsymbol x
\\
\nonumber
+ &\bilinasp(\dispp^n,\bsp) + \bilinafp(\velp^n,\bfp) 
-  \bilinbfp(\pp^n,\bfp) + \bilinbfp(\bpp,\velp^n) + \bilinbsp(\bpp,\dt\dispp^n) -  \bilinbsp(\pp^n,\bsp) 
\\
\nonumber
+& \int_{\interf} \gamma_f \mu_f h^{-1} \dt \dispp^n \cdot \ttau_p\, \bm  \cdot \ttau_p\
+ \int_{\interf} \gamma_f \mu_f h^{-1}  (\velp^n + \aalpha \dt \dispp^n)\cdot \nn_p\,  (\bfp+ \aalpha \bsp) \cdot \nn_p
\\
\nonumber
= & - \int_{\interf} \nn_p \cdot \strf^{n-1} \nn_p\, \big(-\bfp - \aalpha \bsp \big)\cdot \nn_p
- \int_{\interf} \ttau_p \cdot \strf^{n-1} \nn_p \big( - \bsp) \cdot \ttau_p
\\
\label{eq:explicit-solid-b}
+ & \int_{\interf} \gamma_f \mu_f h^{-1} \velf^{n-1} \cdot \ttau_p\, \bsp \cdot \ttau_p
+ \int_{\interf} \gamma_f \mu_f h^{-1} \velf^{n-1} \cdot \nn_p\, \big(\bfp + \aalpha \bsp\big) \cdot \nn_p.
\end{align}

\item[$(a.2)$] Solve the \emph{fluid problem}: using the fluid domain unit normal vectors on $\partial\domf$, namely $\nn_f,\ttau_f$ respectively, given $\velp^n,\pp^n,\dispp^n\in \femvp \times \fempp \times \fems$ from the solution of the problem above, we aim to find $\velf^n,\pf^n  \in \femvf \times \fempf$ such that
\begin{align}
\nonumber
& \rho_f \int_{\domf} \dt \velf^n \cdot \bff d\boldsymbol x
+\bilinaff(\velf^n,\bff) - \bilinbff(\pf^n,\bff) + \bilinbff(\bpf,\velf^n) + \stabffp(\dt\pf,\bpf)
\\
\nonumber
- & \int_{\interf} \nn_f \cdot \strf(\sym\bff,-\bpf) \nn_f\, \big( \velf^n - \velp^n \big)\cdot \nn_f
- \int_{\interf} \ttau_f \cdot \strf(\sym\bff,-\bpf) \nn_f \velf^n \cdot \ttau_f
\\
\nonumber
+ & \int_{\interf} \gamma_f \mu_f h^{-1} \velf^n \cdot \nn_f\, \bff \cdot \nn_f
+ \int_{\interf} \gamma_f \mu_f h^{-1} \velf^n  \cdot \ttau_f\, \bff \cdot \ttau_f
\\
\nonumber
= & \int_{\interf} \nn \cdot \strf^{n-1} \nn_f\, \bff \cdot \nn_f
+ \int_{\interf} \ttau \cdot \strf^{n-1} \nn_f \bff \cdot \ttau_f
\\
\nonumber
+ & \int_{\interf} \nn \cdot \strf(\sym\bff,-\bpf) \nn_f\, (-\velp^n - \aalpha \dt \dispp^n)\cdot \nn_f
+ \int_{\interf} \ttau \cdot \strf(\sym\bff,-\bpf) \nn_f (- \dt\dispp^n ) \cdot \ttau_f
\\
\label{eq:explicit-fluid-b}
+ & \int_{\interf} \gamma_f \mu_f h^{-1} (\velp^n + \aalpha \dt \dispp^n )\cdot \nn_f\, \bff \cdot \nn_f
+ \int_{\interf} \gamma_f \mu_f h^{-1} ( \dt \dispp^n ) \cdot \ttau_f\, \bff \cdot \ttau_f,
+ \fff(t_n;\bff)
\end{align}
where $\stabffp(\dt\pf,\bpf)$ is a stabilization term proposed in \cite{MR2498525} acting on the free fluid pressure, that helps to restore the stability of the explicit time advancing scheme,
\begin{equation*}
\stabffp(\dt\pf,\bpf) := \gamma_{stab} \frac{h \tau}{\gamma_f \mu_f} \int_{\interf} \dt \pf^n\,\bpf\,.
\end{equation*}
We also notice that the terms involving $\strf$ are evaluated at the previous time step, to improve the stability of the explicit coupling.

Looking at the interface terms of equations \eqref{eq:explicit-solid-b}, \eqref{eq:explicit-fluid-b}, we observe that Algorithm (a) corresponds to a \emph{Dirichlet-Neumann} splitting, where Neumann type interface conditions, namely \eqref{CBNSF}-\eqref{CBFT}, are assigned to the structure problem and the kinematic conditions \eqref{CNF}-\eqref{CBJS} are enforced when solving the fluid problem. It has been shown in \cite{causin2005added} that this decomposition may lead to unconditionally unstable schemes. However, as originally highlighted in \cite{MR2498525}, this is not the case when Nitsche's enforcement of the interface conditions is adopted, in combination with appropriate stabilization operators, such as $\stabffp(\dt\pf,\bpf)$.
\end{description}

%>>>>>>>>>>>>>>>>>>>>>>>>>>>>>>>>>>>>>>>>>>>>>>>>>>>>>>>>>>>>>>>>>>>>>>>>>>>>>>>>>>>>>>>>>>>>>>>>>>>>>>>>>>>>>>>>>>>>>>>>>>>>>>>>>

If we finally proceed to solve all the problems independently, we obtain the explicit algorithm reported below. We also formulate explicitly the governing and interface conditions that are enforced in practice when each sub-problem is solved.

\medskip

\noindent{\bf Algorithm ($b$):}
\begin{description}
\item[$(b.1)$] given $\velf^{n-1},\pf^{n-1},\velp^{n-1},\pp^{n-1}$ find $\dispp^n,\dotdispp^n$ in $\domp$ such that
\begin{align*}
& \rho_p \int_{\domp} \dt \dotdispp^n \cdot \bsp
+ \rho_p \int_{\domp} \big( \dotdispp^n - \dt \dispp^n \big) \cdot \dotbsp
+ \bilinasp(\dispp^n,\bsp) 
\\
+& \int_{\interf} \gamma_f \mu_f h^{-1} \dt \dispp^n \cdot \ttau_p\, \bsp  \cdot \ttau_p\
+ \int_{\interf} \gamma_f \mu_f h^{-1}  \aalpha \dt \dispp^n \cdot \nn_p\,  \aalpha \bsp \cdot \nn_p
\\
=& \bilinbsp(\pp^{n-1},\bsp) 
- \int_{\interf} \nn_p \cdot \strf^{n-1} \nn_p\, \big( - \aalpha \bsp \big)\cdot \nn_p
- \int_{\interf} \ttau_p \cdot \strf^{n-1} \nn_p \big( - \bsp) \cdot \ttau_p
\\
+ & \int_{\interf} \gamma_f \mu_f h^{-1} \velf^{n-1} \cdot \ttau_p\, \bsp \cdot \ttau_p
+ \int_{\interf} \gamma_f \mu_f h^{-1} \big( \velf^{n-1} - \velp^{n-1} \big)\cdot \nn_p\, \aalpha \bsp \cdot \nn_p.
\end{align*}
This problem is equivalent to solving the elastodynamics equation, namely \eqref{B1}, where the pressure term has been time-lagged, complemented with the following Robin-type boundary condition on $\Gamma$:
\begin{align*}
\nn_p \cdot \boldsymbol {\sigma}_p \nn_p &= \nn_p \cdot (\boldsymbol{\sigma}_f)^{n-1} \nn_p 
- \gamma_f \mu_f h^{-1} \left(\frac{D \boldsymbol U}{D t} - \boldsymbol{v}^{n-1} + \boldsymbol{q}^{n-1} \right) \cdot \nn_p,
& \text{on} \ \Gamma,
\\
\ttau_p \cdot \boldsymbol {\sigma}_p \nn_p &= \ttau_p \cdot (\boldsymbol{\sigma}_f)^{n-1} \nn_p 
- \gamma_f \mu_f h^{-1} \left(\frac{D \boldsymbol U}{D t} - \boldsymbol{v}^{n-1} \right) \cdot \ttau_p,
& \text{on} \ \Gamma.
\end{align*}

\item[$(b.2)$] given $\velf^{n-1},\pf^{n-1}$ and $\dispp^n$, find $\velp,\pp^n$ in $\domp$ such that
\begin{align*}
& s_0 \int_{\domp} \dt \pp^n \bpp d\boldsymbol x 
+ \bilinafp(\velp^n,\bfp) - \bilinbfp(\pp^n,\bfp) + \bilinbfp(\bpp,\velp^n) 
+ \stabffq\big(\dt\velp  \cdot \nn_p, \bfp\cdot\nn_p \big)
\\
+ & \int_{\interf} \gamma_f \mu_f h^{-1}  \velp^n \cdot \nn_p\, \bfp \cdot \nn_p
=  - \bilinbsp(\bpp,\dt\dispp^n)
\\
+ & \int_{\interf} \gamma_f \mu_f h^{-1} \big( \velf^{n-1} - \aalpha \dt \dispp^{n-1} \big)\cdot \nn_p\, \bfp \cdot \nn_p
+ \int_{\interf} \nn_p \cdot \strf^{n-1} \nn_p\, \bfp \cdot \nn_p.
\end{align*}
This problem consists of the \emph{dual-mixed} weak form of Darcy equations \eqref{B2}-\eqref{B3} complemented with the following boundary condition,
\begin{equation*}
p_p = - \nn_p \cdot (\boldsymbol{\sigma}_f)^{n-1} \nn_p 
- \gamma_f \mu_f h^{-1} \left( \boldsymbol{v}^{n-1} - \aalpha \left(\frac{D \boldsymbol{U}}{D t}\right)^{n-1} \right)\cdot \nn_p,
\quad \text{on} \ \Gamma.
\end{equation*}

\item[$(b.3)$] given $\velp^n,\pp^n,\dispp^n$, find $\velf^n,\pf^n$ in $\domf$ such that
\begin{align*}
& \rho_f \int_{\domf} \dt \velf^n \cdot \bff d\boldsymbol x
+\bilinaff(\velf^n,\bff) - \bilinbff(\pf^n,\bff) + \bilinbff(\bpf,\velf^n) 
\\
+ & \stabffp(\dt\pf,\bpf) + \stabffv\big(\dt\velf^n \cdot \nn_f, \bff\cdot\nn_f \big) 
\\
-& \int_{\interf} \strf(\sym\bff,-\bpf) \nn_f \cdot \velf^n 
+ \int_{\interf} \gamma_f \mu_f h^{-1} \velf^n \cdot \bff 
= \fff(t_n;\bff) + \int_{\interf} \strf^{n-1} \nn_f \cdot \bff 
\\
- & \int_{\interf} \ttau_f \cdot \strf(\sym\bff,-\bpf) \nn_f \dt\dispp^n \cdot \ttau_f
- \int_{\interf} \nn_f \cdot \strf(\sym\bff,-\bpf) \nn_f\, \big(\velp^n + \aalpha \dt \dispp^n \big)\cdot \nn_f
\\
+ & \int_{\interf} \gamma_f \mu_f h^{-1} \big( \velp^n + \aalpha \dt \dispp^n \big)\cdot \nn_f \, \bff \cdot \nn_f
+ \int_{\interf} \gamma_f \mu_f h^{-1} \dt \dispp^n \cdot \ttau_f \, \bff \cdot \ttau_f.
\end{align*}
Before time-lagging of the term $\strf^{n-1} \nn_f \cdot \bff$, this problem corresponds to the fluid equations \eqref{NS1}-\eqref{NS2},
where the kinematic conditions,
\begin{equation}
\boldsymbol{v} \cdot \nn_f = \left(\boldsymbol{q} + \aalpha \frac{D \boldsymbol U}{D t} \right) \cdot \nn_f,
\quad
\boldsymbol v \cdot \boldsymbol \ttau_f = \frac{D\boldsymbol U}{D t} \cdot \boldsymbol \ttau_f 
\quad \textrm{on} \; \Gamma,
\end{equation}
have been enforced using the classical Nitsche's method formulation for boundary conditions \cite{hansbo2005nitsche}.
We observe that new stabilization terms $\stabffv,\stabffq$ have been introduced into the problem formulation. 
Their role is to control the increment of $\velf^n,\velp^n$ over two subsequent time steps, namely
\begin{align*}
\stabffq\big(\dt \velp^n \cdot \nn, \bfp\cdot\nn \big) 
& = \gamma_{stab}^\prime \gamma_{f} \mu_f \frac{\tau}{h} \int_{\interf} \dt\velp \cdot \nn \bfp\cdot\nn\, ,
\\
\stabffv\big(\dt\velf^n \cdot \nn, \bff\cdot\nn \big) 
& = \gamma_{stab}^\prime \gamma_{f} \mu_f \frac{\tau}{h} \int_{\interf} \dt\velf^n \cdot \nn \bff\cdot\nn\, .
\end{align*}
\end{description}

%>>>>>>>>>>>>>>>>>>>>>>>>>>>>>>>>>>>>>>>>>>>>>>>>>>>>>>>>>>>>>>>>>>>>>>>>>>>>>>>>>>>>>>>>>>>>>>>>>>>>>>>>>>>>>>>>>>>>>>>>>>>>>>>>>

%>>>>>>>>>>>>>>>>>>>>>>>>>>>>>>>>>>>>>>>>>>>>>>>>>>>>>>>>>>>>>>>>>>>>>>>>>>>>>>>>>>>>>>>>>>>>>>>>>>>>>>>>>>>>>>>>>>>>>>>>>>>>>>>>>
\subsection{Stability analysis}

Following Burman \& Fernandez, \cite{MR2498525}, the stability analysis of the explicit time advancing scheme is more easily carried out by rewriting the method as a fully implicit scheme. As a result, suitable iteration residuals appear in the equations. Then, the explicit method will share the same stability properties of the implicit one, provided that the residuals can be controlled by means of the energy of the system.

Denoting with $\aasi(\cdot,\cdot)$ the collection of terms on the left hand side of equation \eqref{eq:DIP}, that is the monolithic problem formulation, with $\yy=\{\dispp,\dotdispp,\velf,\pf,\velp,\pp\}$ the vector of all the solution components and with $\zzh$ the corresponding discrete test function, the loosely coupled scheme is equivalent to solve a discrete problem governed by $\aasi(\cdot,\cdot)$ complemented by suitable residual terms. Below, we denote these additional terms by $\TT^n_{*}, *=1,2,3,4$. For the sake of clarity, the dependence of terms $\TT^n_{*}$ on at least one of the basis functions $\bff,\bfp,\bsp$ has been explicitly put into evidence, while the dependence on the numerical solution is hidden, but only the reference time step, $n$, is reported. As a result, the loosely couple scheme is equivalent to find $\yy^n$ such that for any $\zzh$ the following equation is satisfied,
\begin{align}
\label{eq:imex}
& \aasi(\yy^n,\zzh)
+ \stabffp(\dt\pf^n,\bpf) 
+ \vartheta \big( \stabffq\big(\dt \velp^n \cdot \nn, \bfp\cdot\nn \big) 
+ \stabffv\big(\dt\velf^n \cdot \nn, \bff\cdot\nn \big) \big)
\\
\nonumber
& \text{(term denoted as $\TT^n_{1,a}(\bff,\bfp,\bsp)$)}
\\
\nonumber
& = \fff^n(t_n;\zzh) 
+ \int_{\hGamma} \gamma_f \mu_f h^{-1} \big(\velf^{n}- \velf^{n-1}\big) \cdot \ttau\, \big( - \bsp \big)\cdot \ttau
\\
\nonumber
& \text{(terms denoted as $\TT^n_{1,b}(\bff,\bfp,\bsp)$ and $\TT^{n,\vartheta}_{1,c}(\bff,\bfp,\bsp)$)}
\\
\nonumber
& + \int_{\interf} \vartheta \gamma_f \mu_f h^{-1} \big( (\velf^n-\velf^{n-1}) - \aalpha \dt (\dispp^n - \dispp^{n-1}) \big) \cdot \nn (-\bfp) \cdot \nn
\\
\nonumber
& + \int_{\interf} \gamma_f \mu_f h^{-1} \big( (\velf^{n}-\velf^{n-1}) - \vartheta (\velp^{n}-\velp^{n-1}) \big)\cdot \nn\, \big(-(1-\vartheta) \bfp -\aalpha \bsp \big)\cdot \nn
\\
\nonumber
& \text{(terms denoted as $\TT^n_{2}(\bff,\bfp,\bsp)$)}
\\
\nonumber
& + \int_{\interf} 2 \mu_f \nn \cdot \big(\DD(\velf^n) - \DD(\velf^{n-1}) \big) \nn\, \big(\bff - \bfp - \aalpha \bsp \big)\cdot \nn 
\\
\nonumber
& + \int_{\interf} 2 \mu_f \ttau \cdot \big(\DD(\velf^n) - \DD(\velf^{n-1}) \big) \nn (\bff - \bsp ) \cdot  \ttau
\\
\nonumber
& \text{(term denoted as $\TT^n_{3}(\bff,\bfp,\bsp)$ and $\TT^{n,\vartheta}_{4}(\bff,\bfp,\bsp)$)}
\\
\nonumber
& - \int_{\interf} \nn \cdot \big( \strfp^{n} - \strfp^{n-1} \big) \nn\, \big(\bff - \bfp - \aalpha \bsp \big)\cdot \nn
   - \int_{\domp} \alpha \vartheta \big(\pp^{n} - \pp^{n-1}\big) \nabla \cdot \bsp .
\end{align}

If the two step algorithm is addressed, namely Algorithm (a), then the parameter $\vartheta$ in \eqref{eq:imex} is $\vartheta=0$, while it becomes $\vartheta=1$ when the fully uncoupled scheme is considered, that is Algorithm ($b$).  
To prove the stability of the loosely coupled time advancing scheme, we show that an energy estimate similar to \eqref{eq:energy-coupled} holds true also in this case. To this purpose, we apply to equation \eqref{eq:imex} the same test functions used for the proof of Theorem \ref{th:energy}. Then, summing up with respect to the index $n$ and multiplying by $\tau$ and using an argument similar to the one for \eqref{eq:energy-coupled} we obtain:
\begin{align}
\nonumber
& \energyf^N + \energys^N 
+ \tau \frac{\gamma_{stab}}{2} \frac{h}{\gamma_f \mu_f} \|\pf^N\|_{\interf}^2
+ \tau \frac{\gamma_{stab}^\prime}{2} \frac{\gamma_f \mu_f}{h} \vartheta \big( \|\velf^N\cdot\nn_f\|_\interf^2 + \|\velp^N\cdot\nn_p\|_\interf^2 \big)
\\
\nonumber
& + \tau \sum_{n=1}^N \Big[ 2\mu_f \big(1 -\frac{(\sym+1)}{2}\epsfa C_{TI} -\frac{\epsfab}{4} \big) \|\DD(\velf^n)\|^2_{\domf} 
+ \kappa^{-1} \|\velp^n\|^2_{\domp} 
\\
\nonumber
& + \frac{\tau}{2} \Big( \rho_f \|\dt \velf^n\|_{\domf}^2
+ 2 \mu_p \|\dt\DD(\dispp^n)\|^2_{\domp}
+ s_0 \|\dt \pp^n\|_{\domp}^2 
+ \lambda_p \| \dt \nabla \cdot \dispp^n \|_{\domp}^2 \Big)
\\
\nonumber
& +  \big((\gamma_f - (\sym+1)(\epsilon_f^\prime)^{-1}\big) \mu_f h^{-1} \big( \| \big( \velf^n - \velp^n - \aalpha \dt \dispp^n \big)\cdot \nn \|_{\interf}^2
+ \| \big(\velf^n - \dt \dispp^n \big) \cdot \ttau \|_{\interf}^2 \big)
\\
\nonumber
& + \frac{\gamma_{stab}}{2}\frac{ h}{\gamma_f \mu_f} \|\pf^n-\pf^{n-1}\|_{\interf}^2 
  + \frac{\gamma_{stab}^\prime}{2} \frac{\gamma_f \mu_f}{h} \vartheta
\big(\|(\velf^n - \velf^{n-1})\cdot\nn_f\|_\interf^2 + \|(\velp^n - \velp^{n-1})\cdot\nn_p\|_\interf^2 \big)
\Big]
\\
\nonumber
& \leq \energyf^0 + \energys^0 
+ \tau \sum_{n=1}^N \big[ 
\TT^n_{1,a}+\TT^n_{1,b}+\TT^{n,\vartheta}_{1,c}+\TT^n_{2}+\TT^n_{3}+\TT^{n,\vartheta}_{4}
\big](\velf^n,\velp^n,\dt\dispp^n)
\\
\label{eq:stability-1}
& + \frac{\gamma_{stab}}{2} \frac{\tau h}{\gamma_f \mu_f}  \|\pf^0\|_{\interf}^2
+ \frac{\gamma_{stab}^\prime}{2} \gamma_f \mu_f \frac{\tau}{h} \vartheta
\big(\|\velf^0 \cdot\nn_f\|_\interf^2 + \|\velp^0 \cdot\nn_p\|_\interf^2 \big)
+ \tau \sum_{n=1}^N (2\epsfab \mu_f)^{-1} \|\fff(t_n)\|^2,
\end{align}
where, owing to \eqref{eq:dtt}, we have applied the following identities to manipulate the stabilization terms,
\begin{align*}
\tau \sum_{n=1}^N  \stabffp(\dt\pf^n,\pf^n) 
&=\tau \sum_{n=1}^N \frac{\gamma_{stab}}{2} \frac{h}{\gamma_f \mu_f} \|\pf^n-\pf^{n-1}\|_{\interf}^2
\\
&+ \tau \frac{\gamma_{stab}}{2} \frac{h}{\gamma_f \mu_f} \Big( \|\pf^N\|_{\interf}^2 - \|\pf^0\|_{\interf}^2 \Big),
\\
\tau \sum_{n=1}^N \stabffv\big(\dt\velf^n \cdot \nn_f, \velf^n \cdot \nn_f \big) 
&=\tau \sum_{n=1}^N \frac{\gamma_{stab}^\prime}{2} \frac{\gamma_f \mu_f}{h} \|(\velf^n - \velf^{n-1})\cdot\nn_f\|_\interf^2
\\
&+\tau \frac{\gamma_{stab}^\prime}{2} \frac{\gamma_f \mu_f}{h} \Big(\|\velf^N\cdot\nn_f\|_\interf^2 - \|\velf^0\cdot\nn_f\|_\interf^2 \Big) ,
\\
\tau \sum_{n=1}^N \stabffq\big(\dt\velp^n \cdot \nn_p, \velp \cdot \nn_p \big)
&=\tau \sum_{n=1}^N \frac{\gamma_{stab}^\prime}{2} \frac{\gamma_f \mu_f}{h} \|(\velp^n - \velp^{n-1})\cdot\nn_p\|_\interf^2
\\
&+\tau \frac{\gamma_{stab}^\prime}{2} \frac{\gamma_f \mu_f}{h} \Big(\|\velp^N\cdot\nn_p\|_\interf^2 - \|\velp^0\cdot\nn_p\|_\interf^2 \Big).
\end{align*}

The stability of the considered loosely coupled schemes follows from \eqref{eq:stability-1} provided that the terms T1, T2, T3 and T4, defined in equation \eqref{eq:imex}, are bounded by the energy of the discrete problem. To perform this analysis, we extend the approach proposed by  Burman \& Fernandez in \cite{MR2498525} to the case where the flow interacts with a poroelastic structure. As a result, the residual terms defined in \eqref{eq:imex} are more complicated. Their bounds are summarized in the following Lemma.
\begin{lemma}\label{lemma:residuals}
For any $\epsfaa,\epsfaaa,\epsfiv,\epspaa > 0$ the following upper bound holds true,
where $C_{TI}$ is the constant of inequality \eqref{eq:trace_inverse},
\begin{align}
\label{eq:stability-2}
&\tau \sum_{n=1}^N \big[ 
\TT^n_{1,a}+\TT^n_{1,b}+\TT^{n,\vartheta}_{1,c}+\TT^n_{2}+\TT^n_{3}+\TT^{n,\vartheta}_{4}
\big](\velf^n,\velp^n,\dt\dispp^n)
\\ 
\nonumber
&\leq \frac{\gamma_f \mu_f}{2} \frac{\tau}{h} 
\big( (1+\vartheta) \| \velf^{0} \cdot \nn \|_{\interf}^2 + \| \velf^{0} \cdot \ttau \|_{\interf}^2
+ \vartheta \| \velp^{0} \cdot \nn \|_{\interf}^2 + \vartheta \| \aalpha \dt \dispp^{0} \cdot \nn \|_{\interf}^2 \big)
+ \frac{\epsfaa C_{TI}}{2}\mu_f  \|\DD(\velf^0)\|_{\domf}^2  
\\
\nonumber
& + \tau \sum_{n=1}^N \Big[ \frac{\mu_f}{2h} \big(\gamma_f(\epsfaaa + \epsfiv + \vartheta) + (\epsfaa)^{-1} \big) 
\|\big(\velf^{n} - \velp^{n} - \aalpha \dt \dispp^n \big)\cdot \nn\|_{\interf}^2 
+ \frac{\gamma_f \mu_f}{2h}  \|\big(\velf^{n} - \dt \dispp^n \big)\cdot \ttau\|_{\interf}^2 
\\
\nonumber
& + \epsfaa C_{TI}\mu_f  \|\DD(\velf^n)\|_{\domf}^2  
+ \vartheta \Big( \frac{\alpha^2}{s_0} \tau \| \dt \nabla \cdot \dispp^n \|_{\domp}^2
+ \frac{s_0 \tau}{4} \|\dt \pp^n\|_{\domp}^2 \Big)
\\
\nonumber
& + \frac{(\epsfiv)^{-1}}{2} \frac{h}{\gamma_f \mu_f} \|\pf^n-\pf^{n-1}\|_{\interf}^2
+ \frac{\vartheta \big( (\epsfaaa)^{-1} - 1 \big)}{2} \frac{\gamma_f \mu_f}{h}
\Big( \| (\velf^{n}-\velf^{n-1}) \cdot \nn\|_{\interf}^2 + \| (\velp^{n}-\velp^{n-1}) \cdot \nn\|_{\interf}^2 \Big)
\Big].
\end{align}
\end{lemma}

{\bf Proof:}
\begin{description}
\item To estimate $\TT^n_{1,a}$ we proceed as follows:
\begin{multline*}
\TT^n_{1,a}(\velf^n,\velp^n,\dt\dispp^n) 
 = \int_{\interf} \gamma_f \mu_f h^{-1} (\velf^{n}- \velf^{n-1}) \cdot \ttau\, \big(\velf^{n} - \dt \dispp^n \big)\cdot \ttau
- \int_{\interf} \gamma_f \mu_f h^{-1} (\velf^{n}- \velf^{n-1}) \cdot \ttau\,\velf^{n} \cdot \ttau
\\
 \leq \frac{\gamma_f \mu_f}{2} h^{-1} \|\big(\velf^{n} - \dt \dispp^n \big)\cdot \ttau\|_{\interf}^2
+ \frac{\gamma_f \mu_f}{2} \frac{\tau}{h} \Big( \tau \| \dt \velf^{n} \cdot \ttau \|_{\interf}^2 - \dt \| \velf^{n} \cdot \ttau\|_{\interf}^2 - \tau \| \dt \velf^{n} \cdot \ttau\|_{\interf}^2 \Big).
\end{multline*}
When the fully uncoupled scheme is used, the control of the residual $\TT^n_{1,b}$ is analogous to the one of $\TT^n_{1,a}$. 
For simplicity of notation, we introduce $\velz = \velf - \aalpha\dt\dispp$. Then we obtain,
\begin{align*}
\TT^n_{1,b}(\velf^n,\velp^n,\dt\dispp^n)
&= \int_{\interf} \vartheta \gamma_f \mu_f h^{-1} (\velz^{n} - \velz^{n-1}) \cdot \nn\, \big(\velz^{n} - \velp^n \big)\cdot \nn 
- \int_{\interf} \vartheta \gamma_f \mu_f h^{-1} (\velz^{n} - \velz^{n-1}) \cdot \nn \velz^{n} \cdot \nn
\\
& \leq \frac{\vartheta \gamma_f \mu_f}{2} h^{-1} \|\big(\velf^{n} - \velp^{n} - \aalpha \dt \dispp^n \big)\cdot \nn\|_{\interf}^2
\\
& + \frac{\vartheta \gamma_f \mu_f}{2} \frac{\tau}{h} 
\Big( \tau  \| \dt \velz^{n} \cdot \nn \|_{\interf}^2 - \dt \| \velz^{n} \cdot \nn\|_{\interf}^2 - \tau \| \dt \velz^{n} \cdot \nn\|_{\interf}^2 \Big)
\\
& = \frac{\vartheta \gamma_f \mu_f}{2} h^{-1} \|\big(\velf^{n} - \velp^{n} - \aalpha \dt \dispp^n \big)\cdot \nn\|_{\interf}^2
- \frac{\vartheta \gamma_f \mu_f}{2} \frac{\tau}{h} \dt \| (\velf^{n}-\aalpha\dt\dispp^{n}) \cdot \nn\|_{\interf}^2.
\end{align*}
To estimate $\TT^{n,\vartheta}_{1,c}$ we introduce the auxiliary variable $\velr = \velf - \vartheta \velp$ at the interface $\Gamma$. 
Then, $\TT^{n,\vartheta}_{1,c}$ can be rewritten as:
\begin{align*}
\TT^{n,\vartheta}_{1,c}(\velf^n,\velp^n,\dt\dispp^n) 
&= \int_{\interf} \gamma_f \mu_f h^{-1} (\velr^{n} - \velr^{n-1}) \cdot \nn\, \big(\velr^{n} - (1-\vartheta) \velp^n - \aalpha \dt \dispp^n \big)\cdot \nn 
\\
& - \int_{\interf} \gamma_f \mu_f h^{-1} (\velr^{n} - \velr^{n-1}) \cdot \nn \velr^{n} \cdot \nn
\\
& \leq \epsfaaa \frac{\gamma_f \mu_f}{2} h^{-1} \|\big(\velf^{n} - \velp^{n} - \aalpha \dt \dispp^n \big)\cdot \nn\|_{\interf}^2
\\
& + \frac{\gamma_f \mu_f}{2} \frac{\tau}{h} \Big( \tau (\epsfaaa)^{-1}  \| \dt \velr^{n} \cdot \nn \|_{\interf}^2 - \dt \| \velr^{n} \cdot \nn\|_{\interf}^2 - \tau \| \dt \velr^{n} \cdot \nn\|_{\interf}^2 \Big)
\\
& \leq  \frac{\gamma_f \mu_f}{2} h^{-1} \Big(
\epsfaaa\|\big(\velf^{n} - \velp^{n} - \aalpha \dt \dispp^n \big)\cdot \nn\|_{\interf}^2
\\
& + \big( (\epsfaaa)^{-1} - 1 \big) \| (\velf^{n}-\velf^{n-1}) \cdot \nn\|_{\interf}^2
  + \vartheta \big( (\epsfaaa)^{-1} - 1 \big) \| (\velp^{n}-\velp^{n-1}) \cdot \nn\|_{\interf}^2
\Big)
\\
& - \frac{\gamma_f \mu_f}{2} \frac{\tau}{h} \dt \| (\velf^{n}- \vartheta \velp^{n}) \cdot \nn\|_{\interf}^2.
\end{align*}
As result of that 
\begin{align*}
&\tau \sum_{n=1}^N \TT^n_{1}(\velf^n,\velp^n,\dt\dispp^n) 
= \tau \sum_{n=1}^N \big[ \TT^n_{1,a} + \TT^n_{1,b} + \TT^{n,\vartheta}_{1,c} \big](\velf^n,\velp^n,\dt\dispp^n) 
\\
&\leq \frac{\gamma_f \mu_f}{2} \frac{\tau}{h} 
\big( (1+\vartheta) \| \velf^{0} \cdot \nn \|_{\interf}^2 + \| \velf^{0} \cdot \ttau \|_{\interf}^2
+ \vartheta \| \velp^{0} \cdot \nn \|_{\interf}^2 + \vartheta \| \aalpha \dt \dispp^{0} \cdot \nn \|_{\interf}^2 \big)
\\
&+ \tau \sum_{n=1}^N \Big[ \frac{\gamma_f \mu_f}{2h} 
\Big( \|\big(\velf^{n} - \dt \dispp^n \big)\cdot \ttau\|_{\interf}^2
+ (\vartheta+\epsfaaa)\|\big(\velf^{n} - \velp^{n} - \aalpha \dt \dispp^n \big)\cdot \nn\|_{\interf}^2 \Big)
\\
&+ \big( (\epsfaaa)^{-1} - 1 \big) \| (\velf^{n}-\velf^{n-1}) \cdot \nn\|_{\interf}^2
+ \vartheta \big( (\epsfaaa)^{-1} - 1 \big) \| (\velp^{n}-\velp^{n-1}) \cdot \nn\|_{\interf}^2.
\Big].
\end{align*}

\item For terms $\TT^n_{2}$ and $\TT^n_{3}$ we proceed as in the derivation of the energy estimate:
\begin{align*}
\TT^n_{2}(\velf^n,\velp^n,\dt\dispp^n) 
& \leq \epsfaa  C_{TI} \mu_f \big( \|\DD(\velf^n)\|_{{\Omega}_f}^2 + \|\DD(\velf^{n-1})\|_{\domf}^2 \big)
\\
& + (\epsfaa)^{-1} h^{-1} \mu_f \big( \|\big(\velf^{n} - \velp^{n} - \aalpha \dt \dispp^n \big)\cdot \nn\|_{\interf}^2
+ \|\big(\velf^{n} - \dt \dispp^n \big)\cdot \ttau\|_{\interf}^2 \big),
\\[13pt]
\TT^n_{3}(\velf^n,\velp^n,\dt\dispp^n)
& \leq (\epsfiv)^{-1} \frac{h}{2\gamma_f \mu_f} \|\pf^{n} - \pf^{n-1}\|_{\interf}^2 
+ \epsfiv \frac{\gamma_f \mu_f}{2} h^{-1} \|\big(\velf^{n} - \velp^{n} - \aalpha \dt \dispp^n \big)\cdot \nn\|_{\interf}^2.
\end{align*}

\item To determine an appropriate upper bound for $\TT^{n,\vartheta}_{4}$ it is important to take into account of the magnitude of the physical parameters $s_0,\,\alpha,\,\lambda_p$, for which we have approximately $s_0 \ll 1 \simeq \alpha \ll \lambda_p$. Then, for the first of terms $\TT^{n,\vartheta}_{4}$ we have:
\begin{equation*}
\TT^{n,\vartheta}_{4}(\velf^n,\velp^n,\dt\dispp^n) 
= \vartheta \alpha \tau \int_{\domp} \big( \dt \pp \big)\, \big( \dt \nabla \cdot \dispp^n \big) 
\leq \vartheta \frac{s_0 }{2} \tau \|\dt \pp^n\|_{\domp}^2 + \vartheta \frac{\alpha^2}{2s_0} \tau \| \dt \nabla \cdot \dispp^n \|_{\domp}^2.
\end{equation*}
\hfill$\Box$
\end{description}

The stability of the loosely coupled scheme is just a direct consequence of \eqref{eq:stability-1} combined with \eqref{eq:stability-2}, 
which gives
\begin{align}
\nonumber
& \energyf^N + \energys^N 
+ \tau \frac{\gamma_{stab}}{2} \frac{h}{\gamma_f \mu_f} \|\pf^N\|_{\interf}^2
+ \tau \vartheta \frac{\gamma_{stab}^\prime}{2} \frac{\gamma_f \mu_f}{h} \big( \|\velf^N\cdot\nn_f\|_\interf^2 + \|\velp^N\cdot\nn_p\|_\interf^2 \big)
\\
\nonumber
& + \tau \sum_{n=1}^N \Big[ \mu_f (2 -(\sym+1)\epsfa C_{TI} -\frac{\epsfab}{2}  -2\epsfaa C_{TI}) \|\DD(\velf^n)\|^2_{\domf} 
+ \kappa^{-1} \|\velp^n\|^2_{\domp} 
\\
\nonumber
& + \frac{\tau}{2} \Big( \rho_f \|\dt \velf^n\|_{\domf}^2
+ 2 \mu_p \|\dt\DD(\dispp^n)\|^2_{\domp}
+ (1-\vartheta) s_0 \|\dt \pp^n\|_{\domp}^2 
+ \lambda_p \big( 1 - \vartheta \frac{\alpha^2}{\lambda_p s_0} \big) \| \dt \nabla \cdot \dispp^n \|_{\domp}^2 \Big)
\\
\nonumber
& +  \frac12 \big(\gamma_f(2-(\epsfaaa+\epsfiv+\vartheta)) - 2(\sym+1)(\epsfa)^{-1} - 2(\epsfaa)^{-1}\big) \mu_f h^{-1} 
\big( \| \big( \velf^n - \velp^n - \aalpha \dt \dispp^n \big)\cdot \nn \|_{\interf}^2
\\
\nonumber
& + \frac12 \big(\gamma_f - 2(\sym+1)(\epsfa)^{-1} - 2(\epsfaa)^{-1}\big) \| \big(\velf^n - \dt \dispp^n \big) \cdot \ttau \|_{\interf}^2 \big)
+ \frac12 (\gamma_{stab}-(\epsfiv)^{-1})\frac{ h}{\gamma_f \mu_f} \|\pf^n-\pf^{n-1}\|_{\interf}^2 
\\
\nonumber
& + \frac12 (\gamma_{stab}^\prime-((\epsfaaa)^{-1}-1)) \frac{\gamma_f \mu_f}{h} 
\big(\|(\velf^n - \velf^{n-1})\cdot\nn_f\|_\interf^2 + \vartheta \|(\velp^n - \velp^{n-1})\cdot\nn_p\|_\interf^2 \big)
\Big]
\\
\nonumber
& \leq \energyf^0 + \energys^0 
+ \frac{\epsfaa C_{TI}}{2} \mu_f  \|\DD(\velf^0)\|_{\domf}^2  
\\
\nonumber
& + \frac{\gamma_f \mu_f}{2} \frac{\tau}{h} 
\big( (1+\vartheta) \| \velf^{0} \cdot \nn \|_{\interf}^2 + \| \velf^{0} \cdot \ttau \|_{\interf}^2
+ \vartheta \| \velp^{0} \cdot \nn \|_{\interf}^2 + \vartheta \| \aalpha \dt \dispp^{0} \cdot \nn \|_{\interf}^2 \big)
\\
\label{eq:stability-3}
& + \vartheta \frac{\gamma_{stab}^\prime}{2} \gamma_f \mu_f \frac{\tau}{h}
\big(\|\velf^0 \cdot\nn_f\|_\interf^2 + \|\velp^0 \cdot\nn_p\|_\interf^2 \big)
+ \frac{\gamma_{stab}}{2} \frac{\tau h}{\gamma_f \mu_f}  \|\pf^0\|_{\interf}^2
+ \tau \sum_{n=1}^N (2\epsfab \mu_f)^{-1} \|\fff(t_n)\|^2.
\end{align}

\null

The stability of the loosely coupled scheme is summarized by the following result.
\begin{theorem}\label{th:stability}
For any $\epsfa,\epsfab,\epsfaa,\epsfaaa,\epsfiv>0$ such that
\begin{equation*}
(2 -(\sym+1)\epsfa C_{TI}  -2\epsfaa C_{TI} -\epsfab/2) > 0,
\quad 
(2-(\epsfaaa+\epsfiv+\vartheta)) = \delta > 0,
\end{equation*}
provided that the penalty and stabilization parameters are large enough, more precisely
\begin{gather*}
\gamma_{stab} \geq (\epsfiv)^{-1}, 
\quad 
\gamma_{stab}^\prime \geq (\epsfaaa)^{-1} - 1,
\\
\gamma_f > \frac{2(\sym+1)(\epsfa)^{-1} + 2(\epsfaa)^{-1}}{\delta},
\end{gather*}
that the problem coefficients satisfy the following restriction
\begin{equation}\label{eq:stab_constraint}
\vartheta \frac{\alpha^2}{\lambda_p s_0} < 1,
\end{equation}
and that, asymptotically, the ratio $\tau / h$ is bounded, namely that the CFL-type condition $\tau < C h$ is satisfied, for a given positive constant $C$, there exist constants $0<c_f^\prime < 1< C_f^\prime$, uniformly independent of the mesh characteristic size, such that
\begin{align}
\nonumber
& \energyf^N + \energys^N
+ \tau \frac{\gamma_{stab}}{2} \frac{h}{\gamma_f \mu_f} \|\pf^N\|_{\interf}^2
+ \tau \vartheta \frac{\gamma_{stab}^\prime}{2} \frac{\gamma_f \mu_f}{h} \big( \|\velf^N\cdot\nn_f\|_\interf^2 + \|\velp^N\cdot\nn_p\|_\interf^2 \big)
\\
\nonumber
&+  c_f^\prime \tau \sum_{n=1}^N \Big[ 
\mu_f \|\DD(\velf^n)\|^2_{L^2(\domf)} 
+ \kappa^{-1} \|\velp^n\|^2_{L^2(\domp)} 
+ \mu_f h^{-1} \big( \| ( \velf^n - \velp^n - \aalpha \dt \dispp^n )\cdot \nn \|_{\interf}^2
+ \| \big(\velf^n - \dt \dispp^n \big) \cdot \ttau \|_{\interf}^2 \big)\
\\
\nonumber
& + \tau \big( \rho_f \|\dt \velf^n\|_{\domf}^2 
+ 2 \mu_p \|\dt\DD(\dispp^n)\|^2_{\domp} 
+ s_0 \|\dt \pp^n\|_{\domp}^2 
+ \lambda_p \| \dt \nabla \cdot \dispp^n \|_{\domp}^2
\big) \Big]
\\
\nonumber
& \leq \energyf^0 + \energys^0
+ C_f^\prime \mu_f^{-1} \tau \sum_{n=1}^N \|\fff(t_n)\|^2
\\
\label{eq:stability-4}
& + C_f^\prime \mu_f \big( \| \velf^{0} \cdot \nn \|_{\interf}^2 + \| \velf^{0} \cdot \ttau \|_{\interf}^2
+ \| \velp^{0} \cdot \nn \|_{\interf}^2 + \| \aalpha \dt \dispp^{0} \cdot \nn \|_{\interf}^2 
+ \|\DD(\velf^0)\|^2_{L^2(\domf)} + \mu_f^{-2} \|\pf^0\|_{\interf}^2\big).
\end{align}
More precisely, we have
\begin{align*}
c_f^\prime &< \min \{ 
(2 - (\sym+1) \epsfa C_{TI} -2\epsfaa C_{TI}-\epsfab/2),
\big( 1 - \vartheta \frac{\alpha^2}{\lambda_p s_0} \big),
\ldots
\\
&\ldots\big(\gamma_f(2-(\epsfaaa+\epsfiv+\vartheta)) - 2(\sym+1)(\epsfa)^{-1} - 2(\epsfaa)^{-1}\big),
\big(\gamma_f - 2(\sym+1)(\epsilon_f^\prime)^{-1} - 2(\epsfaa)^{-1}\big)\},
\\
C_f^\prime &> \max \{(2\epsfab)^{-1}, \frac{\epsfaa C_{TI}}{2}, \frac{\gamma_{stab}}{2} \frac{\tau h}{\gamma_f}, \vartheta \frac{\gamma_{stab}^\prime}{2} \gamma_f \}.
\end{align*}
\end{theorem}
\begin{corollary}\label{th:stability-1}
When the Biot problem is not decoupled, namely when $\vartheta=0$, the stability constraints can be relaxed as follows:
\begin{itemize}
\item[-] inequality \eqref{eq:stab_constraint} is always satisfied.
\item[-] it is possible to choose $\epsfaaa=1$ and $0< \epsfiv <1$ to satisfy $(2-(\epsfaaa+\epsfiv+\vartheta) = \delta > 0$. 
As a result $\gamma_{stab}^\prime=0$ is an admissible choice, namely the stabilization operators $\stabffq$ and $\stabffv$ are no longer needed.
\end{itemize}
\end{corollary}
The numerical experiments that will be presented later on suggest that the loosely couple scheme in Algorithm $(b)$ turns out to be stable in practice even when the stabilization terms $\stabffq$ and $\stabffv$ are omitted, namely when we choose $\gamma_{stab}^\prime = 0$. This behavior can be explained by the following stability result, which is a variant of Theorem \ref{th:stability}.
\begin{corollary}\label{th:stability-2}
Under the ranges of $\epsfa,\epsfab,\epsfaa,\epsfaaa,\epsfiv,\gamma_f,\gamma_{stab}$ that satisfy the assumptions of Theorem \ref{th:stability} and provided that \eqref{eq:stab_constraint} holds true, the following additional conditions are sufficient to guarantee stability when $\vartheta=1$ and $\gamma_{stab}^\prime=0$:
\begin{gather*}
\exists C>0 \ s.t. \ \tau < C h^2 \ and \ \rho_f > \frac12 ((\epsfaaa)^{-1}-1)) \gamma_f \mu_f C_{TI} C ,
\\
\kappa^{-1} > ((\epsfaaa)^{-1}-1)) C_{TI} h^{-2} .
\end{gather*}
\end{corollary}
{\bf Proof:} the result follows by isolating the last term of \eqref{eq:stability-2}. It provides an upper bound in terms of the energy nom, without resorting to the stabilization operators $\stabffq$ and $\stabffv$. In particular, we have
\begin{equation*}
h^{-1} \| (\velf^{n}-\velf^{n-1}) \cdot \nn\|_{\interf}^2
=  \frac{\tau^2}{h} \| \dt \velf^{n} \cdot \nn\|_{\interf}^2
\leq \frac{\tau^2}{h^2} C_{TI} \| \dt \velf^{n} \|_{\domf}^2
\leq C C_{TI} \tau \| \dt \velf^{n} \|_{\domf}^2,
\end{equation*}
provided that $\tau < C h^2$.  As a result we obtain,
\begin{equation*}
\frac{\gamma_f \mu_f}{h} \frac{\vartheta \big( (\epsfaaa)^{-1} - 1 \big)}{2} \| (\velf^{n}-\velf^{n-1}) \cdot \nn\|_{\interf}^2
\leq \gamma_f \mu_f C C_{TI} \frac{\vartheta \big( (\epsfaaa)^{-1} - 1 \big)}{2} \tau \| \dt \velf^{n} \|_{\domf}^2.
\end{equation*}
The term above can be combined with $\tau \rho_f \|\dt \velf^n\|_{\domf}^2$ on the left hand side of \eqref{eq:stability-4}.
For the second part of the considered residuals we have,
\begin{equation*}
\sumtau h^{-1} \| (\velp^{n}-\velp^{n-1}) \cdot \nn\|_{\interf}^2 \leq \sumtau 2 C_{TI} h^{-2} \|\velp^{n}\|_{\domp}^2 + \frac{\tau}{h} \|\velp^{0}\|_{\domp}^2,
\end{equation*}
that is 
\begin{equation*}
\frac{\vartheta \big( (\epsfaaa)^{-1} - 1 \big)}{2} \frac{\gamma_f \mu_f}{h} \sumtau \| (\velp^{n}-\velp^{n-1}) \cdot \nn\|_{\interf}^2 
\leq \gamma_f \mu_f \frac{\vartheta \big( (\epsfaaa)^{-1} - 1 \big)}{2} \Big( \sumtau 2 C_{TI} h^{-2} \|\velp^{n}\|_{\domp}^2 + \frac{\tau}{h} \|\velp^{0}\|_{\domp}^2 \Big),
\end{equation*}
which can be in turn combined with $\kappa^{-1} \|\velp^n\|^2_{\domp}$. \hfill $\Box$

%>>>>>>>>>>>>>>>>>>>>>>>>>>>>>>>>>>>>>>>>>>>>>>>>>>>>>>>>>>>>>>>>>>>>>>>>>>>>>>>>>>>>>>>>>>>>>>>>>>>>>>>>>>>>>>>>>>>>>>>>>>>>>>>>>

%>>>>>>>>>>>>>>>>>>>>>>>>>>>>>>>>>>>>>>>>>>>>>>>>>>>>>>>>>>>>>>>>>>>>>>>>>>>>>>>>>>>>>>>>>>>>>>>>>>>>>>>>>>>>>>>>>>>>>>>>>>>>>>>>>
\subsection{Splitting error analysis}

Adopting a loosely coupled scheme may affect the accuracy of the underlying approximation methods, in particular for time discretization. In this section we highlight these additional sources of error and analyse how they affect the accuracy of the scheme. To keep the technical details of the analysis to a minimum, without disregarding a rigorous approach, we mainly focus on the time discretization error. As a result, the analysis consists of comparing the following problems:
\begin{align*}
\text{Semidiscrete problem (SDP):} & \text{ find } \yy(t) \ s.t. \ \aas(\yy(t),\zzh) = \fff(t;\zzh)
\\
\text{Discrete problem, implicit coupling (DIP):} & \text{ find } \yyin \ s.t. \ \aasi(\yyin,\zzh) = \fff^n(\zzh)
\\
\text{Discrete problem, explicit coupling (DEP):} & \text{ find } \yyen \ s.t. \ \aase(\yyen,\zzh) = \fff^n(\zzh)
\end{align*}
where for simplicity of notation we denote with $\yy=\{\dispp,\dotdispp,\velf,\pf,\velp,\pp\}$ the vector of all the solution components and with $\zzh$ the corresponding discrete test function. At the same time, we introduce a compact notation for the norms that we have used to prove the stability of (DIP) and (DEP), in Theorems \ref{th:energy} and \ref{th:stability}, respectively. For the one used in Theorem \ref{th:energy} we have
\begin{align*}
\triple{\yy}_{\heartsuit,\,N}^2 &:= \energyf^N(\yy) + \energys^N(\yy)
\\
\triple{\yy}_{\nabla,\,n}^2 &:=
2 \mu_f \|\DD(\velf^n)\|^2_{\domf} + \kappa^{-1} \|\velp^n\|^2_{\domp} 
\\
&+ \frac{\tau}{2} \Big( \rho_f \|\dt \velf^n\|_{\domf}^2 
+ 2 \mu_p \|\dt \DD(\dispp^n)\|^2_{\domp} 
+ s_0 \|\dt \pp^n\|_{\domp}^2
+ \lambda_p \| \dt \nabla \cdot \dispp^n \|_{\domp}^2 \Big)
\\
&+ \mu_f h^{-1} \big( \| \big( \velf^n - \velp^n - \aalpha \dt \dispp^n \big)\cdot \nn \|_{\interf}^2
+ \| \big(\velf^n - \dt \dispp^n \big) \cdot \ttau \|_{\interf}^2 \big).
\end{align*}
We proceed similarly for the norm of Theorem \ref{th:stability}, for which $\triple{\cdot}_{\nabla,\,n}$ remains unchanged and we define
\begin{equation*}
\triple{\yy}_{\Diamond,\,N}^2:= \triple{\yy}_{\heartsuit,\,N}^2
+ \frac{\tau h}{\gamma_f \mu_f} \|\pf^N\|_{\interf}^2
+ \vartheta \gamma_f \mu_f \frac{\tau}{h} \big( \|\velf^N\cdot\nn_f\|_\interf^2 + \|\velp^N\cdot\nn_p\|_\interf^2 \big).
\end{equation*}

Our analysis stems form the assumption that the \emph{discrete problem with fully implicit coupling is asymptotically consistent with respect to time discretization}. In particular, let problem (SDP) be rearranged such that
\begin{equation}\label{eq:residual-i}
\aas(\yy(t_n),\zzh) = \aasi(\yy(t_n),\zzh) - \rsi(\yy(t_n),\zzh) = \fff(t_n;\zzh).
\end{equation}
As a result, subtracting (DIP) from (SDP) we obtain the following equation for $\ein := \yy(t_n) - \yyin$,
\begin{equation*}
\aasi(\ein,\zzh) = \big[ \fff(t_n;\zzh) - \fff^n(\zzh) \big] + \rsi(\yy(t_n),\zzh).
\end{equation*}
\begin{assumption}\label{dip:consistency}
There exists $C_{DIP}>0$ such that, for any $\epsilon > 0$ the following estimates hold true
\begin{equation*}
|\fff(t;\ein) - \fff^n(\ein)| \leq C_{DIP} \sumtau \big( \frac{1}{2\epsilon} \triple{\fff(t_n;\ein)-\fff^n(\ein)}^2 + \frac{\epsilon}{2} \triple{\ein}_{\nabla,\,n}^2 \big) ,
\end{equation*}
\begin{equation*}
\rsi(\yy(t_n),\ein) \leq C_{DIP} \sumtau \big( \frac{1}{2\epsilon} \triple{\rsi(\yy(t_n))}^2 + \frac{\epsilon}{2} \triple{\ein}_{\nabla,\,n}^2 \big) .
\end{equation*}
Furthermore, there exists $p>0$ such that
\begin{equation*}
\sumtau \triple{\fff(t_n;\ein)-\fff^n(\ein)}^2 = {\cal O}(\tau^{2p})
\quad \text{and} \quad
\sumtau \triple{\rsi(\yy(t_n))}^2 = {\cal O}(\tau^{2p}).
\end{equation*}
\end{assumption}
The previous properties combined with Theorem \ref{th:energy} ensure that the time discretization error of (DIP) converges with order $p$ in the norm $\triple{\cdot}_{\heartsuit,\,N}+\sumtau \triple{\cdot}_{\nabla,\,n}$. We also assume that the trace, Poincar\'e and Korn inequalities, respectively reported below,
\begin{equation*}
%\label{trace}
\|\velf\|_{\Gamma}
\leq C_T  \|\velf\|_{\domf}^{1/2}  \| \nabla \velf \|_{\domf}^{1/2},
\quad
%\label{PF}
\|\velf\|_{\domf} \leq C_P \| \nabla \velf\|_{\domf},
\quad
%\label{korn}
\| \nabla \velf \|_{\domf} \leq C_K \| D(\velf) \|_{\domf},
\end{equation*}
where constants $C_P, C_T$ and $C_K$ depend on the domain $\Omega$, see \cite{MR2050138} for details, are satisfied by the finite element approximation space for the free flow velocity.

Starting from this background, we now analyze the time discretization error of (DEP). We observe that the bilinear form corresponding to (DEP) can be split into the following parts,
\begin{equation}\label{eq:splitting-e}
\aase(\yyen,\zzh) = \aasi(\yyen,\zzh) + \stab(\yyen,\zzh) + \spliterr(\yyen,\zzh),
\end{equation}
where $\stab(\cdot,\cdot)$ denotes the collection of the previously defined stabilization terms 
\[\stab(\cdot,\cdot)=s_{f,p}(\cdot,\cdot)+\vartheta \big( s_{f,v}(\cdot,\cdot)+s_{f,q}(\cdot,\cdot) \big)\]
and $\spliterr(\cdot,\cdot)$ is the splitting error residual, as defined in \eqref{eq:imex}.
Furthermore, proceeding as in \eqref{eq:residual-i}, we rearrange (SDP) as follows,
\begin{equation}\label{eq:residual-e}
\aas(\yy(t_n),\zzh) = \aase(\yy(t_n),\zzh) - \rse(\yy(t_n),\zzh) = \fff(t_n;\zzh),
\end{equation}
where $\rse(\cdot,\cdot)$ is the time discretization residual relative to (DEP).
Combining \eqref{eq:splitting-e} and \eqref{eq:residual-e} with (SDP) we obtain
\begin{equation*}
\aasi(\yy(t_n),\zzh) +  \stab(\yy(t_n),\zzh) + \spliterr(\yy(t_n),\zzh) = \fff(t_n;\zzh) + \rse(\yy(t_n),\zzh).
\end{equation*}
which, owing to \eqref{eq:residual-i}, is equivalent to the following additive decomposition of $\rse(\cdot,\cdot)$
\begin{equation*}
\rse(\yy(t_n),\zzh) = \rsi(\yy(t_n),\zzh) +  \stab(\yy(t_n),\zzh) + \spliterr(\yy(t_n),\zzh).
\end{equation*}
In conclusion, subtracting (DEP) from (SDP), we obtain the following equation that combines the time discretization error relative to (DEP), namely $\een := \yy(t_n) - \yyen$, with the residuals quantifying the consistency error due to stabilization terms, equation splitting and time discretization respectively,
\begin{equation}\label{eq:error-e}
\aase(\een,\zzh) = \big[ \fff(t_n;\zzh) - \fff^n(\zzh) \big]
+ \stab(\yy(t_n),\zzh) + \spliterr(\yy(t_n),\zzh) + \rsi(\yy(t_n),\zzh).
\end{equation}
As a consequence of the stability estimate for the explicit scheme, we know that $\yyen$ converges to $\yy(t_n)$, provided that the right hand side of \eqref{eq:error-e} vanishes when $\tau,h \rightarrow 0$.

The asymptotic good behavior of $\rsi(\yy(t_n),\zzh)$ has been previously discussed, as a consequence of the assumption that time discretization of (DIP) is consistent. In what follows, we focus our attention on the analysis of the residuals related to the loosely coupled scheme, namely $\stab(\yy(t_n),\zzh)$ and $\spliterr(\yy(t_n),\zzh)$. 
The different components of the time discretization error are denoted as 
$\een:=\{\edispp^n,\edotdispp^n,\evelf^n,\epf^n,\evelp^n,\epp^n\}$.
The results are summarized in the following lemmas.

\begin{lemma}\label{lemma:stabilization}
For any $\epsilon> 0$, the stabilization terms $s_{f,p}(\cdot,\cdot),\,s_{f,v}(\cdot,\cdot),\,s_{f,q}(\cdot,\cdot)$ satisfy the following upper bounds, respectively,
\begin{equation*}
\sumtau s_{f,p}(\dt \pf(t_n),\epf^n) \leq \frac{\gamma_{stab}}{2} \frac{h\tau}{\gamma_f\mu_f} 
\Big[ 
\epsilon^{-1} \int_\Gamma \int_0^T |\partial_t \pf(t)|^2 dt\,ds + \epsilon \sumtau \|\epf^n\|_\Gamma^2
\Big],
\end{equation*}
\begin{equation*}
\sumtau s_{f,v}(\dt \velf(t_n),\evelf^n) \leq \frac{\gamma_{stab}^\prime}{2} \gamma_f \mu_f \frac{\tau}{h} 
\Big[ 
\epsilon^{-1} \int_\Gamma \int_0^T |\partial_t \velf(t)\cdot\nn|^2 dt\,ds + \epsilon \sumtau \|\evelf^n\cdot\nn\|_\Gamma^2
\Big],
\end{equation*}
\begin{equation*}
\sumtau s_{f,q}(\dt \velp(t_n),\evelp^n) \leq \frac{\gamma_{stab}^\prime}{2} \gamma_f \mu_f \frac{\tau}{h} 
\Big[ 
\epsilon^{-1} \int_\Gamma \int_0^T |\partial_t \velp(t)\cdot\nn|^2 dt\,ds + \epsilon \sumtau \|\evelp^n\cdot\nn\|_\Gamma^2
\Big].
\end{equation*}
As a result, there exists a constant $C_S = \max \{ \gamma_{stab}, \gamma_{stab}^\prime \}$ such that
\begin{equation*}
 \sumtau \stab(\yy(t_n),\een) \leq C_S \sumtau \big( \epsilon^{-1} \triple{\stab(\yy(t_n)}^2 + \epsilon \triple{\een}_{\Diamond,\,n}^2 \big)
\end{equation*}
\begin{multline*}
\sumtau \triple{\stab(\yy(t_n)}^2 \leq \frac{h\tau}{\gamma_f\mu_f} \int_\Gamma \int_0^T |\partial_t \pf(t)|^2 dt\,ds
\\
+ \vartheta \gamma_f \mu_f \frac{\tau}{h} \int_\Gamma \int_0^T \big( |\partial_t \velf(t)\cdot\nn|^2 + |\partial_t \velp(t)\cdot\nn|^2 \big) dt\,ds .
\end{multline*}
\end{lemma}

The splitting error residuals have already been analyzed in Lemma \ref{lemma:residuals} to prove the stability of (DEP). There is however a substantial difference between Lemma \ref{lemma:residuals} and the upper bounds needed to ensure convergence. Here, the splitting error residual $\spliterr(\yy(t_n),\zzh)$ must be tested using $\zzh=\een$. As a consequence, part of the estimates derived in Lemma \ref{lemma:residuals} are no longer valid. Following a different approach, we are able to prove the result below.

\begin{lemma}\label{lemma:spliterr}
For any $\epsilon>0$ the splitting error residuals satisfy the following upper bounds
\begin{align*}
& \sumtau \spliterr(\yy(t_n),\een)  
\\
& \leq \frac{\epsilon}{2} \gamma_f \mu_f \sumtau \Big[
3 h^{-1} \big( \|(\evelf^n -\evelp^n -\aalpha\dt\edispp^n)\cdot\nn\|_\Gamma^2 + \|(\evelf^n - \dt\edispp^n)\cdot\ttau\|_\Gamma^2 \big)
\\
&+ 2 C_T^2 C_P C_K^2 \frac{\tau}{h} \|\DD(\evelf^n)\|_{\domf}^2 
+ \theta \frac{\tau}{h} \big(  \|\dt \edispp^n\cdot\nn\|_\Gamma^2 
+ \|\evelp^n\cdot\nn\|_\Gamma^2 \big) \Big]
\\
& + \frac{(1+\tau)}{2\epsilon} \gamma_f \mu_f \frac{\tau}{h} \Big[ \int_\Gamma \int_0^T 
\big( (1+\vartheta) |\partial_t \velf(t)\cdot\nn|^2 + |\partial_t \velf(t) \cdot \ttau|^2 + \vartheta |\partial_t \velp(t)\cdot\nn|^2 \big) dt\,ds
\\
& + \vartheta\sumtau \|\partial^2_{tt} \dispp(t_{n-1})\cdot\nn\|_\Gamma^2  \Big]
\\
& + \frac{\mu_f \tau^2}{\epsilon \gamma_f} C_{TI} \int_{\domf} \int_0^T |\partial_t \DD(\velf(t))|^2 dt\,d\boldsymbol x
+ \frac{h\tau^2}{2 \epsilon \gamma_f \mu_f} \int_\Gamma \int_0^T |\partial_t \pf(t)|^2 dt\,ds
\\
& + \vartheta \frac{\tau s_0}{2} \int_{\domp} \int_0^T |\partial_t \pp(t)|^2 dt\,d\boldsymbol x
+ \vartheta \frac{\tau \alpha^2}{2 s_0} \int_{\domp} \int_0^T |\partial_t \nabla \cdot \dispp(t)|^2 dt\,d\boldsymbol x .
\end{align*}
As a result, under the restriction of the time discretization step $\tau \leq C h$, there exists a constant $C_R$ such that
\begin{equation*}
\sumtau \spliterr(\yy(t_n),\een)  \leq C_R \sumtau \Big[
\epsilon \triple{\een}_{\nabla,\,n}^2
+\epsilon \triple{\een}_{\Diamond,\,n}^2
+\epsilon^{-1} \triple{\spliterr(\yy(t_n)}^2 \Big]
\end{equation*}
\begin{align*}
&\sumtau \triple{\spliterr(\yy(t_n)}^2 \leq 
\gamma_f \mu_f \frac{\tau}{h} \Big[ 
\vartheta\sumtau \|\partial^2_{tt} \dispp(t_{n-1})\cdot\nn\|_\Gamma^2
\\
&\qquad + \int_\Gamma \int_0^T 
\big( (1+\vartheta) |\partial_t \velf(t)\cdot\nn|^2 + |\partial_t \velf(t) \cdot \ttau|^2 + \vartheta |\partial_t \velp(t)\cdot\nn|^2 \big) dt\,ds
 \Big]
\\
+ & \tau^2 \Big[ \frac{\mu_f }{\gamma_f} C_{TI} \int_{\domf} \int_0^T |\partial_t \DD(\velf(t))|^2 dt\,d\boldsymbol x
+ \frac{h}{\gamma_f \mu_f} \int_\Gamma \int_0^T |\partial_t \pf(t)|^2 dt\,ds \Big]
\\
+ & \vartheta \tau \Big[ s_0\int_{\domp} \int_0^T |\partial_t \pp(t)|^2 dt\,d\boldsymbol x
+ \frac{\alpha^2}{s_0} \int_{\domp} \int_0^T |\partial_t \nabla \cdot \dispp(t)|^2 dt\,d\boldsymbol x \Big].
\end{align*}
\end{lemma}

The proof of Lemma \ref{lemma:spliterr} is postponed at the end of this section.
From the inequality above we observe that the splitting error residuals generate three types of consistency errors, respectively proportional to
\[
\frac{\tau}{h}, \quad \tau^2, \quad \vartheta \tau.
\]
It is clear that the first and the last items of the list are sub-optimally convergent when $\tau,h \rightarrow 0$. In particular, the first item is not even asymptotically convergent in general. Such errors arise when splitting the penalty terms introduced for enforcing the kinematic coupling conditions by Nitsche's method. A similar behavior has been observed in \cite{NME:NME4607}. An additional restriction on $\tau,h$ will be introduced later on to override this drawback. 

\smallskip

To sum up, choosing $\epsilon$ small enough in the upper bounds of the splitting error residuals, analyzed in Lemmas \ref{lemma:stabilization}, \ref{lemma:spliterr}, and combining them with the stability of (DEP), reported in Theorem \ref{th:stability}, we can easily prove that there exist positive constants $c_E=\max\{C_{DIP},C_R,C_S\}$  and $C_E=\epsilon^{-1} \max\{C_{DIP},C_R,C_S\}$  such that
\begin{multline*}
\triple{\eeN}_{\Diamond,N}^2 + \sumtau \triple{\een}_{\nabla,n}^2 
\leq c_E \tau \sum_{n=1}^{N-1} \triple{\een}_{\Diamond,n}^2
\\
+ C_E \sumtau \Big[ 
   \triple{\fff(t_n)-\fff^n}^2
+ \triple{\rsi(\yy(t_n))}^2
+ \triple{\stab(\yy(t_n)}^2
+ \triple{\spliterr(\yy(t_n)}^2
\Big].
\end{multline*}
As a consequence of the discrete Gronwall Lemma we obtain the following error estimate. 
\begin{theorem}\label{th:convergence}
Under the assumptions of Theorem \ref{th:stability} and Lemma \ref{lemma:spliterr}, the following result holds true
\begin{multline*}
\triple{\eeN}_{\Diamond,N}^2 + \sumtau \triple{\een}_{\nabla,n}^2 
\\
\leq C_E \exp{T} \sumtau \Big[ \triple{\fff(t_n)-\fff^n}^2 + \triple{\rsi(\yy(t_n))}^2 + \triple{\stab(\yy(t_n)}^2 + \triple{\spliterr(\yy(t_n)}^2 \Big].
\end{multline*}
If the assumption \ref{dip:consistency} is satisfied together with the restriction of the time discretization step size $\tau \leq c h^2$, then the right hand side is such that
\begin{equation*}
\sumtau \Big[ \triple{\fff(t_n)-\fff^n}^2 + \triple{\rsi(\yy(t_n))}^2 + \triple{\stab(\yy(t_n)}^2 + \triple{\spliterr(\yy(t_n)}^2 \Big] = {\cal O} (\tau^\frac12).
\end{equation*}
\end{theorem}

We observe that the constraint $\tau \leq c h^2$ arises to make sure that the residuals
\[
\frac{\tau}{h}\int_\Gamma \int_0^T 
\big( (1+\vartheta) |\partial_t \velf(t)\cdot\nn|^2 + |\partial_t \velf(t) \cdot \ttau|^2 + \vartheta |\partial_t \velp(t)\cdot\nn|^2 \big) dt\,ds
\]
are asymptotically convergent when $\tau,h \rightarrow 0$. 
Furthermore, these terms converge with a sub-optimal rate ${\cal O} (\tau^\frac12)$.
This drawback is common to the two variants of loosely coupled schemes that we have considered, indeed the terms above do not vanish if we set $\vartheta=0$. This behavior completely agrees with what was observed in \cite{NME:NME4607}. In addition, the following residuals that scale with a suboptimal rate with respect to $\tau,h$ appear when elasticity and Darcy equations are split in the Biot problem,
\[
\vartheta \tau \Big[ s_0\int_{\domp} \int_0^T |\partial_t \pp(t)|^2 dt\,d\boldsymbol x
+ \frac{\alpha^2}{s_0} \int_{\domp} \int_0^T |\partial_t \nabla \cdot \dispp(t)|^2 dt\,d\boldsymbol x \Big].
\]

{\bf Proof of Lemma \ref{lemma:spliterr}:}
We choose $\zzh=\een$ in \eqref{eq:error-e} and we proceed to analyze each term of the decomposition of $\spliterr$ adopted in Lemma \ref{lemma:residuals}, that is
\[
\sumtau \spliterr(\yy(t_n),\een) = \sumtau  \big( \TT_1 + \TT_2 + \TT_3 + \TT_4 \big)(\yy(t_n);\een),
\]
\begin{multline*}
\sumtau \TT_{1,a}(\yy(t_n);\een) 
\\
\leq \frac12 \gamma_f \mu_f \frac{\tau}{h} \sum_{n=1}^N 
\big[ \epsilon^{-1} (1 + \tau ^{-1}) \|(\velf(t_n)-\velf(t_{n-1}))\cdot\ttau\|_\Gamma^2 
+ \epsilon \big( \|(\evelf^n - \dt\edispp)\cdot\ttau\|_\Gamma^2 
+ \tau \|\evelf^n\cdot\ttau\|_\Gamma^2 \big) \big]
\\
\leq \frac12 \gamma_f \mu_f \frac{\tau}{h} \sum_{n=1}^N 
\big[ \epsilon^{-1} (1 + \tau^{-1}) \tau \int_\Gamma \int_0^T |\partial_t \velf(t)\cdot\ttau|^2 dt\,d\boldsymbol x 
+ \epsilon \big( \|(\evelf^n - \dt\edispp)\cdot\ttau\|_\Gamma^2 
+ \tau C_T^2 C_P C_K^2 \|\DD(\evelf^n)\|_{\domf}^2 \big) \big].
\end{multline*}
Since the semidiscrete solution $\yy(t)$ is time-continuous, terms of type $\|g(t_n)-g(t_{n-1})\|_\Gamma^2$ are bounded as follows
\begin{equation*}
\sum_{n=1}^N \|g(t_n)-g(t_{n-1})\|_\Gamma^2 \leq \tau \int_\Gamma \int_0^T |\partial_t g(t)|^2 dt\,d\boldsymbol x,
\end{equation*}
\begin{multline*}
\sumtau \TT_{1,b}(\yy(t_n);\een) 
\leq \frac{\vartheta}{2} \gamma_f \mu_f \frac{\tau}{h} \sum_{n=1}^N 
\Big[ \epsilon^{-1} (1+\tau^{-1}) \Big( \|(\velf(t_n)-\velf(t_{n-1}))\cdot\nn\|_\Gamma^2 
+ \aalpha\|\dt(\dispp(t_n)-\dispp(t_{n-1}))\cdot\nn\|_\Gamma^2 \Big)
\\
+ \epsilon \tau \|\dt\edispp^n\cdot\nn\|_\Gamma^2 \Big]
\leq \frac{\vartheta}{2\epsilon} \gamma_f \mu_f (1+\tau) \frac{\tau}{h}
\Big[ \int_\Gamma \int_0^T |\partial_t \velf(t)\cdot\nn|^2 dt\,d\boldsymbol x
+ \sumtau \|\partial^2_{tt} \dispp(t_{n-1})\cdot\nn\|_\Gamma^2 \Big]
\\
+ \frac{\vartheta \epsilon}{2} \gamma_f \mu_f \frac{\tau}{h} \sumtau \|\dt\edispp^n\cdot\nn\|_\Gamma^2,
\end{multline*}
\begin{multline*}
\sumtau \TT_{1,c}^\vartheta(\yy(t_n);\een) 
\leq \frac12 \gamma_f \mu_f \frac{\tau}{h} \sum_{n=1}^N 
\Big[ \epsilon^{-1} (1 + \tau^{-1}) \Big( \|(\velf(t_n)-\velf(t_{n-1}))\cdot\nn\|_\Gamma^2 + \vartheta \|(\velp(t_n)-\velp(t_{n-1}))\cdot\nn\|_\Gamma^2 \Big)
\\
+ \epsilon \big( \|(\evelf^n - \evelp^n - \aalpha\dt\edispp)\cdot\nn\|_\Gamma^2
+ \tau \|\evelf^n\cdot\nn\|_\Gamma^2 + \tau \vartheta \|\evelp^n\cdot\nn\|_\Gamma^2 \big) \Big]
\\
\leq \frac12 \gamma_f \mu_f \frac{\tau}{h} \sum_{n=1}^N 
\Big[ \epsilon^{-1} (1 + \tau^{-1}) \tau \Big( 
\int_\Gamma \int_0^T \big( |\partial_t \velf(t)\cdot\nn|^2 + \vartheta |\partial_t \velp(t)\cdot\nn|^2 \Big) dt\,d\boldsymbol x 
\\
+ C_T^2 C_P C_K^2 \frac{\tau}{h} \|\DD(\evelf^n)\|_{\domf}^2
+ \epsilon \big( \|(\evelf^n - \evelp^n - \aalpha\dt\edispp)\cdot\nn\|_\Gamma^2
+ \tau \vartheta \|\evelp^n\cdot\nn\|_\Gamma^2 \big) \Big],
\end{multline*}
\begin{multline*}
\sumtau \TT_2(\yy(t_n);\een) 
\leq \sumtau
\big[ \frac{\mu_f h}{\epsilon \gamma_f} \|\DD(\velf(t_n) - \velf(t_{n-1}))\|_\Gamma^2
\\
+ \frac{\epsilon \gamma_f \mu_f}{h} \big( \|(\evelf^n - \evelp^n - \aalpha\dt\edispp)\cdot\nn\|_\Gamma^2 + \|(\evelf^n - \dt\edispp)\cdot\ttau\|_\Gamma^2 \big]
\leq \frac{\mu_f \tau^2}{\epsilon \gamma_f} C_{TI} \int_{\domf} \int_0^T |\partial_t \DD(\velf(t))|^2 dt\,d\boldsymbol x
\\
+ \epsilon \gamma_f \mu_f\frac{\tau}{h} \sum_{n=1}^N \big( \|(\evelf^n - \evelp^n - \aalpha\dt\edispp)\cdot\nn\|_\Gamma^2 + \|(\evelf^n - \dt\edispp)\cdot\ttau\|_\Gamma^2 \big),
\end{multline*}
\begin{multline*}
\sumtau \TT_3(\yy(t_n);\een) 
\leq \frac{h\tau^2}{2 \epsilon \gamma_f \mu_f} \int_\Gamma \int_0^T |\partial_t \pf(t)|^2 dt\,ds
+\frac{\epsilon}{2} \gamma_f \mu_f \frac{\tau}{h} \sum_{n=1}^N  \|(\evelf^n - \evelp^n - \aalpha\dt\edispp)\cdot\nn\|_\Gamma^2 ,
\end{multline*}
\begin{equation*}
\sumtau \TT_4(\yy(t_n);\een) 
\leq 
\vartheta \frac{\tau s_0}{2} \int_{\domp} \int_0^T |\partial_t \pp(t)|^2 dt\,d\boldsymbol x
+ \vartheta \frac{\tau \alpha^2}{2 s_0} \int_{\domp} \int_0^T |\partial_t \nabla \cdot \dispp(t)|^2 dt\,d\boldsymbol x.
\end{equation*}
\hfill$\Box$

%>>>>>>>>>>>>>>>>>>>>>>>>>>>>>>>>>>>>>>>>>>>>>>>>>>>>>>>>>>>>>>>>>>>>>>>>>>>>>>>>>>>>>>>>>>>>>>>>>>>>>>>>>>>>>>>>>>>>>>>>>>>>>>>>>

%>>>>>>>>>>>>>>>>>>>>>>>>>>>>>>>>>>>>>>>>>>>>>>>>>>>>>>>>>>>>>>>>>>>>>>>>>>>>>>>>>>>>>>>>>>>>>>>>>>>>>>>>>>>>>>>>>>>>>>>>>>>>>>>>>
\subsection{Using the loosely coupled schemes as preconditioners}

The error analysis shows that poor accuracy is the main drawback of the loosely coupled methods (DEP), in contrast to their significant advantage in terms of computational efficiency, if compared to the fully implicit or monolithic scheme (DIP).  This issue has been investigated in \cite{MR2498525,NME:NME4607} for a similar scheme, where the authors propose to use subiterations between the subproblems defining the loosely coupled algorithm, among other options. The numerical investigation addressed in \cite{NME:NME4607} confirms that this approach may cure the aforementioned drawbacks. However, a thorough stability analysis of these more accurate variants of the loosely coupled schemes is not available yet. Furthermore, the results presented in \cite{NME:NME4607} refer to the coupling of a linear flow model with an elastic structure. The extension to poroelastic problem is a completely open subject. 

Another approach to blend the computational efficiency of loosely coupled schemes with the accuracy of the monolithic ones is to use the former schemes as preconditioners for the solution of the latter. As a result, we will end up with a method which shares the desirable properties of the two approaches. As we will show later on, the stability analysis that we have performed for the DIP and DEP schemes suggests that the loosely coupled method may perform as an optimal preconditioner for the solution of the algebraic problem related to DIP. The key factor to prove this property is the spectral equivalence of the stiffness matrices related to DIP and DEP respectively.

%-------------------------------------------------------------------------------------------------------------------------------------------------------------------------------------------------------------

The stability analysis preformed in Theorems \ref{th:energy} and \ref{th:stability} informs us about the spectra of the matrices related to the bilinear forms $\aasi(\cdot,\cdot)$ and $\aase(\cdot,\cdot)$. To start with, let us consider the monolithic scheme complemented with the pressure stabilization operator $\stabffp(\cdot,\cdot)$. Let us denote with $\aasis(\cdot,\cdot)$ the resulting scheme. Considering this variant of DIP will significantly facilitate the forthcoming analysis, but does not affect the stability and convergence properties of the method, as confirmed by the error analysis. Let us denote with $\maasis$ and $\maase$, the generalized stiffness matrices corresponding to the bilinear forms $\aasis$ and $\aase$, respectively. We notice that $\maasis$ and $\maase$ are neither symmetric nor positive definite. Furthermore, let $\YY$ be the vector of degrees of freedom corresponding to the finite element function $\yy$. For the sake of clarity, the block structure of the algebraic monolithic problem is illustrated below:
{\scriptsize
\begin{multline*}
\begin{bmatrix}
\big(M_f + A_f + \Gamma_{f}^\gamma + \Gamma_{f}^\sigma + \sym(\Gamma_{f}^\sigma)^T\big) & -(B_{pf}+\Gamma_{pf})^T &  \Gamma_{qf}^T & 0 & \Gamma_{sf}^T & 0
\\[10pt]
(B_{pf}+\Gamma_{pf}) & S_{p} & \Gamma_{qp}^T & 0 & \Gamma_{sp}^T & 0
\\[10pt]
\Gamma_{qf} & \Gamma_{qp} & (A_q + \Gamma_{q}) & -B_{pq}^T & \Gamma_{sq}^T & 0
\\[10pt]
0 & 0 & B_{pq} & M_p & B_{sp}^T & 0
\\[10pt]
\Gamma_{sf} & \Gamma_{sp}  & \Gamma_{sq} & -B_{sp} & (A_s + \Gamma_{s}) & M_s
\\[10pt]
0 & 0 & 0 & 0 & - M_s & \dot{M}_s
\end{bmatrix}
\cdot
\begin{bmatrix}
\boldsymbol v^n \\[10pt] p_f^n \\[10pt] \boldsymbol q^n \\[10pt] p_p^n \\[10pt] \boldsymbol U^n \\[10pt] \dot{\boldsymbol U}^n
\end{bmatrix}
\\
=
\begin{bmatrix}
\fff(t_n) \\[10pt] 0 \\[10pt] 0 \\[10pt] 0 \\[10pt] 0 \\[10pt] 0 
\end{bmatrix}
+
\begin{bmatrix}
(M_f + S_f) & 0 &  0 & 0 & (\Gamma_{sf}^\gamma+\sym\Gamma_{sf}^\sigma)^T & 0
\\[10pt]
0 & S_{p} & 0 & 0 & (\sym\Gamma_{sp}^\sigma)^T & 0
\\[10pt]
0 & 0 & S_{q} & 0 & (\Gamma_{sq}^\gamma)^T & 0
\\[10pt]
0 & 0 & 0 & M_p & B_{sp}^T & 0
\\[10pt]
0 & 0 & 0 & 0 & \Gamma_{s} & M_s
\\[10pt]
0 & 0 & 0 & 0 & - M_s & 0
\end{bmatrix}
\cdot
\begin{bmatrix}
\boldsymbol v^{n-1} \\[10pt] p_f^{n-1} \\[10pt] \boldsymbol q^{n-1} \\[10pt] p_p^{n-1} \\[10pt] \boldsymbol U^{n-1} \\[10pt] \dot{\boldsymbol U}^{n-1}
\end{bmatrix}.
\end{multline*}
}
We notice that the matrix above has a remarkable structure, combining symmetric and skew symmetric blocks. In particular, all the matrices related to diagonal blocks are symmetric. The loosely coupled scheme is equivalent to the following upper block-triangular system, according to the fact that each sub-problem can be solved independently, but to ensure stability they must be addressed in a precise order,
{\scriptsize
\begin{multline*}
\begin{bmatrix}
\big(M_f + A_f + \Gamma_{f}^\gamma + \sym(\Gamma_{f}^\sigma)^T+ S_{f}\big) & -B_{pf}^T &  \Gamma_{qf}^T & 0 & \Gamma_{sf}^T & 0
\\[10pt]
(B_{pf} + \Gamma_{pf}) & S_{p} & \Gamma_{qp}^T & 0 & \Gamma_{sp}^T & 0
\\[10pt]
0 & 0 & (A_q + \Gamma_{q} + S_{q}) & - B_{pq}^T & \Gamma_{sq}^T & 0
\\[10pt]
0 & 0 & B_{pq} & M_p & B_{sp}^T & 0
\\[10pt]
0 & 0 & 0 & 0 & (A_s + \Gamma_{s}) & M_s
\\[10pt]
0 & 0 & 0 & 0 & - M_s & \dot{M}_s
\end{bmatrix}
\cdot
\begin{bmatrix}
\boldsymbol v^n \\[10pt] p_f^n \\[10pt] \boldsymbol q^n \\[10pt] p_p^n \\[10pt] \boldsymbol U^n \\[10pt] \dot{\boldsymbol U}^n
\end{bmatrix}
\\
=
\begin{bmatrix}
\fff(t_n) \\[10pt] 0 \\[10pt] 0 \\[10pt] 0 \\[10pt] 0 \\[10pt] 0 
\end{bmatrix}
+
\begin{bmatrix}
\big(M_f - \Gamma_{f}^\sigma + S_{f}\big) & \Gamma_{pf}^T &  0 & 0 & \Gamma_{sf}^T & 0
\\[10pt]
0 & S_{p} & 0 & 0 & \Gamma_{sp}^T & 0
\\[10pt]
\Gamma_{qf} & \Gamma_{qp} & S_{q} & 0 & \Gamma_{sq}^T & 0
\\[10pt]
0 & 0 & 0 & M_p & B_{sp}^T & 0
\\[10pt]
\Gamma_{sf} & \Gamma_{sp}  & \Gamma_{sq} & B_{sp} & \Gamma_{s} & M_s
\\[10pt]
0 & 0 & 0 & 0 & - M_s & 0
\end{bmatrix}
\cdot
\begin{bmatrix}
\boldsymbol v^{n-1} \\[10pt] p_f^{n-1} \\[10pt] \boldsymbol q^{n-1} \\[10pt] p_p^{n-1} \\[10pt] \boldsymbol U^{n-1} \\[10pt] \dot{\boldsymbol U}^{n-1}
\end{bmatrix}.
\end{multline*}
}
Before proceeding, we introduce the energy norm related to the bilinear form $\aasis$, which will be used to prove the spectral equivalence of $\aasis$ and $\aase$:
\[
\triple{\yy}^2 := \triple{\yy}_{\heartsuit,\,n}^2 + \tau \triple{\yy}_{\nabla,\,n}^2 + \frac{\tau h}{\gamma_f \mu_f} \|\pf\|_{\interf}^2.
\]
More precisely the term by term definition of the energy norm is the following:
\begin{multline*}
\triple{\yy}^2 := \rho_f \|\velf\|_{\domf}^2 
+ 2 \mu_p \|\DD(\dispp)\|^2_{\domp} 
+ s_0 \|\pp\|_{\domp}^2
+ \lambda_p \| \nabla \cdot \dispp \|_{\domp}^2 
+ \rho_p\| \dotdispp \|^2_{\domp}   
\\
+ 2 \tau \mu_f \|\DD(\velf)\|^2_{\domf} + \tau \kappa^{-1} \|\velp\|^2_{\domp} 
+ \mu_f \frac{\tau}{h} \big( \| \big( \velf - \velp - \aalpha \dt \dispp \big)\cdot \nn \|_{\interf}^2
+ \| \big(\velf - \dt \dispp \big) \cdot \ttau \|_{\interf}^2 \big)
+ \frac{\tau h}{\gamma_f \mu_f} \|\pf\|_{\interf}^2.
\end{multline*}

The stability analysis of the schemes DIP and DEP can be reformulated in order to show that the following properties hold true.

\begin{lemma}\label{l:equivalence}
Under the assumptions of Theorem \ref{th:energy}, $\aasis$ satisfies the following property:
there exist constants $\widetilde{C}>\widetilde{c}>0$ independent on $\tau,h$ such that 
\begin{equation*}
\widetilde{c} \triple{\yy}^2 \leq \aasis(\yy,\yy) \leq \widetilde{C} \triple{\yy}^2 \quad \text{for any} \ \yy\, .
\end{equation*}
Under the assumptions of Theorem \ref{th:stability}, Corollary \ref{th:stability-1} if $\vartheta=0$ or Corollary \ref{th:stability-2} when $\vartheta=1$, $\aase$ satisfies the following property: there exist constants $\widehat{C}>\widehat{c}>0$ independent on $\tau,h$ such that 
\begin{equation*}
\widehat{c} \triple{\yy}^2 \leq \aase(\yy,\yy) \leq \widehat{C} \triple{\yy}^2 \quad \text{for any} \ \yy\, .
\end{equation*}
\end{lemma}

{\bf Proof:} The positivity of $\aasis$ and $\aase$, namely the lower bounds, is a consequence of Theorems \ref{th:energy} and \ref{th:stability}, respectively. The desired results can be shown by applying Theorems \ref{th:energy} and \ref{th:stability} to the particular case of a single time step, using $\fff(t^n;\bff)=0$ and $\yy^{n-1}=0$. In particular, the analysis leading to \eqref{eq:stability-DIP}and \eqref{eq:stability-4} should be adapted using the test functions,
\begin{equation*}
\bff=\velf^n,\,\bfp=\velp^n,\,\bpf=\pf^n,\,\bpp=\pp^n,\,\bsp=\dispp^n,\, \dotbsp=\dotdispp^n.
\end{equation*}

Using the same strategy, the upper bound of $\aasis$ follows naturally form the definition of the bilinear form.
According to equation \eqref{eq:imex}, the upper bound of $\aase$ is a consequence of the previous property combined with inequality \eqref{eq:stability-2}, which in turn must be considered using $\yy^{n-1}=0$. \hfill $\Box$

Lemma \ref{l:equivalence} allows us to prove the following sequence of inequalities,
\begin{equation}\label{eq:equivalence}
\frac{\widetilde{c}}{\widehat{C}} \aase(\yy,\yy) \leq \aasis(\yy,\yy) \leq \frac{\widetilde{C}}{\widehat{c}} \aase(\yy,\yy) \quad \text{for any} \ \yy\, .
\end{equation}
Applying the equivalence between the functional and algebraic representation of the discrete problem,
\[
\aasis(\yy,\yy) = \YY^T \maasis \YY \ \text{and} \ \aase(\yy,\yy) = \YY^T \maase \YY \quad \text{for any} \ \YY,
\]
the inequalities \eqref{eq:equivalence} can be rewritten as follows,
\begin{equation}\label{eq:spectral}
\frac{\widetilde{c}}{\widehat{C}}  \leq \frac{\YY^T \maasis \YY}{\YY^T \maase \YY}  \leq \frac{\widetilde{C}}{\widehat{c}} \quad \text{for any} \ \YY \ \text{s.t.} \ \YY^T \maase \YY \neq 0\, .
\end{equation}
Property \eqref{eq:spectral} can be interpreted as the spectral equivalence of matrices $\maasis$ and $\maase$, and suggests that $\maase$ is a good candidate for preconditioning the monolithic system $\maasis$. This evidence is further supported by the following well known estimate of the GMRES algorithm \cite{eisenstat:1983}, when applied to the generic left-preconditioned system $P^{-1} A \YY = P^{-1} \FF$.
\begin{theorem}\label{t:gmres}
Provided that $P$ is symmetric positive definite, the residuals computed by GMRES at the $k$-th iteration satisfy
\begin{equation}\label{eq:convrate}
\frac{\|P^{-1}\RR^k\|_P}{\|P^{-1}\RR^0\|_P} \leq \big(1-\mu_1\mu_2 \big)^{k/2},
\end{equation}
where $\mu_1$ and $\mu_2$ are defined by
\[
\mu_1 := \min\limits_{\ZZ \neq 0} \frac{\ZZ^T A \ZZ}{\ZZ^T P \ZZ},
\quad
\mu_2 := \min\limits_{\ZZ \neq 0} \frac{\ZZ^T A^{-1} \ZZ}{\ZZ^T P^{-1} \ZZ},
\]
where $\|\cdot\|_P$ denotes the norm induced by $P$.
\end{theorem}
We observe that the constant $\mu_1$ can be heuristically related to $\frac{\widetilde{c}}{\widehat{C}}$ and $\mu_2$ to $\frac{\widehat{c}}{\widetilde{C}}$. Unfortunately, we can not combine this Theorem with \eqref{eq:spectral} to determine a rigorous estimate of the GMRES convergence rate, when it is applied to $(\maase)^{-1} \maasis$. There are two shortcomings of the available properties of the proposed computational scheme: $(i)$ $\maase$ is not symmetric; $(ii)$ $\maase$ is only positive semidefinite. The latter conclusion is a consequence of Lemma \ref{l:equivalence}. It shows that $\triple{\yy}^2$ is equivalent to $\YY^T \maase \YY$. However, $\triple{\yy}^2$ is only a seminorm because it does not provide control on $\|\pf^n\|_{\domf}$, but only on $\|\pf^n\|_{\Gamma}$. In conclusion, \eqref{eq:spectral} only qualitatively suggests that $\maase$ may be a good preconditioner for $\maasis$, and we will thoroughly investigate this conjecture with numerical experiments addressed in the forthcoming sections.

%>>>>>>>>>>>>>>>>>>>>>>>>>>>>>>>>>>>>>>>>>>>>>>>>>>>>>>>>>>>>>>>>>>>>>>>>>>>>>>>>>>>>>>>>>>>>>>>>>>>>>>>>>>>>>>>>>>>>>>>>>>>>>>>>>

\section{Numerical results and discussion}\label{sec:numerics}

The numerical experiments that will be addressed in this section have multiple purposes. On one hand, we aim at supporting the theoretical results on the accuracy and computational efficiency of the proposed methods. On the other hand, we aim to show that the proposed computational framework is suitable for a broad spectrum of applications. In particular, we will focus our attention on prototype problems arising from cardiovascular mechanics and geomechanics. For this reason, this section is mainly split in two parts, each referring to one specific application. All the numerical tests are implemented using the finite element library \textit{FreeFem++} \cite{freefem_new,freefem}.

%>>>>>>>>>>>>>>>>>>>>>>>>>>>>>>>>>>>>>>>>>>>>>>>>>>>>>>>>>>>>>>>>>>>>>>>>>>>>>>>>>>>>>>>>>>>>>>>>>>>>>>>>>>>>>>>>>>>>>>>>>>>>>>>>>
\subsection{Analysis of the scheme applied to fluid-structure interaction in arteries}

We consider a classical benchmark problem used for FSI problems, see \cite{bukavc2012fluid,MR2498525} and references therein, which consists in studying the propagation of a single pressure wave whose amplitude is comparable to the pressure difference between systolic and diastolic phases of the heartbeat. 

\paragraph{Description of the problem.} For simplicity, we adopt a two-dimensional (2D) geometrical model that consists of a 2D elastic structure superposed to a 2D fluid channel. The fluid-structure interaction is excited by the following time-dependent pressure profile, prescribed at the inflow of the channel
\begin{equation*}%\label{pressure}
 p_{in}(t) = \left\{\begin{array}{l@{\ } l} 
\frac{p_{max}}{2} (1 - \cos(\frac{2 \pi t}{T_{max}})) & \; \textrm{if} \; t \le T_{max}, \\
0   & \; \textrm{if} \; t>T_{max},
 \end{array} \right.
\end{equation*}
where $p_{max} = 13334 \; dyne/cm^2$ and $T_{max}=0.003 \; s$. 
To make this test case represent more closely the behavior of an artery, we also slightly modify the governing equation for elastic skeleton as follows,
\[
\rho_{p} \frac{D^2 \boldsymbol U}{D t^2} + \xi \boldsymbol U- \nabla \cdot \boldsymbol \sigma_p = 0.
\]
The additional term $\xi \boldsymbol U$ comes from the axially symmetric formulation, accounting for the recoil due to the circumferential strain. Namely, it acts like a spring term, keeping the top and bottom structure displacements connected in 2D, see, e.g.,~\cite{ma1992numerical}. 

The reference values of the parameters used in this study fall within the range of physiological values for blood flow and are reported in Table~\ref{tab:1}. The propagation of the pressure wave is analyzed over the time interval $[0,0.006]$ s. The final time is selected such that the pressure wave barely reaches the outflow section. In this way, the non-physical reflected waves that will originate at the outflow section for longer simulation times do not pollute the considered results.

\begin{table}[ht!]
\begin{center}
{\small{
\begin{tabular}{l c c c}
\hline
\textbf{Parameter} & \textbf{Symbol} & \textbf{Units} &\textbf{Reference value} \\
\hline
\hline
\textbf{Radius}&  $R$ & (cm)  & $0.5$  \\
\textbf{Length}& $L$ & (cm) & $6$  \\
\textbf{Poroelastic wall thickness}& $r_p$ & (cm) & $0.1$  \\
\textbf{Poroelastic wall density}& $\rho_{p}$ & (g/cm$^3$) & $1.1$  \\
\textbf{Fluid density}& $\rho_f$ & (g/cm$^3$)& $1$\\
\textbf{Dyn. viscosity}& $\mu$ & (g/cm s) & $0.035$    \\
\textbf{Lam\'e coeff.}& $\mu_p $ & (dyne/cm$^2$) & $1.07 \times 10^6$  \\
\textbf{Lam\'e coeff.}& $\lambda_p $ & (dyne/cm$^2$) & $4.28 \times 10^6$ \\
\textbf{Hydraulic conductivity}& $\kappa $ & (cm$^3$ s/g) & $5 \times 10^{-9}$ \\
\textbf{Mass storativity coeff.}& $s_0 $ & (cm$^2$/dyne) & $5 \times 10^{-6}$   \\
\textbf{Biot-Willis constant}& $\alpha$ &  & 1 \\
\textbf{Spring coeff.}& $\xi$ & (dyne/cm$^4$) & $5 \times 10^7$   \\\hline
\end{tabular}
}}
\end{center}
\caption{Geometry, fluid and structure parameters.}
\label{tab:1}
\end{table}

\begin{figure}
\begin{center}
\includegraphics[width=0.85\textwidth]{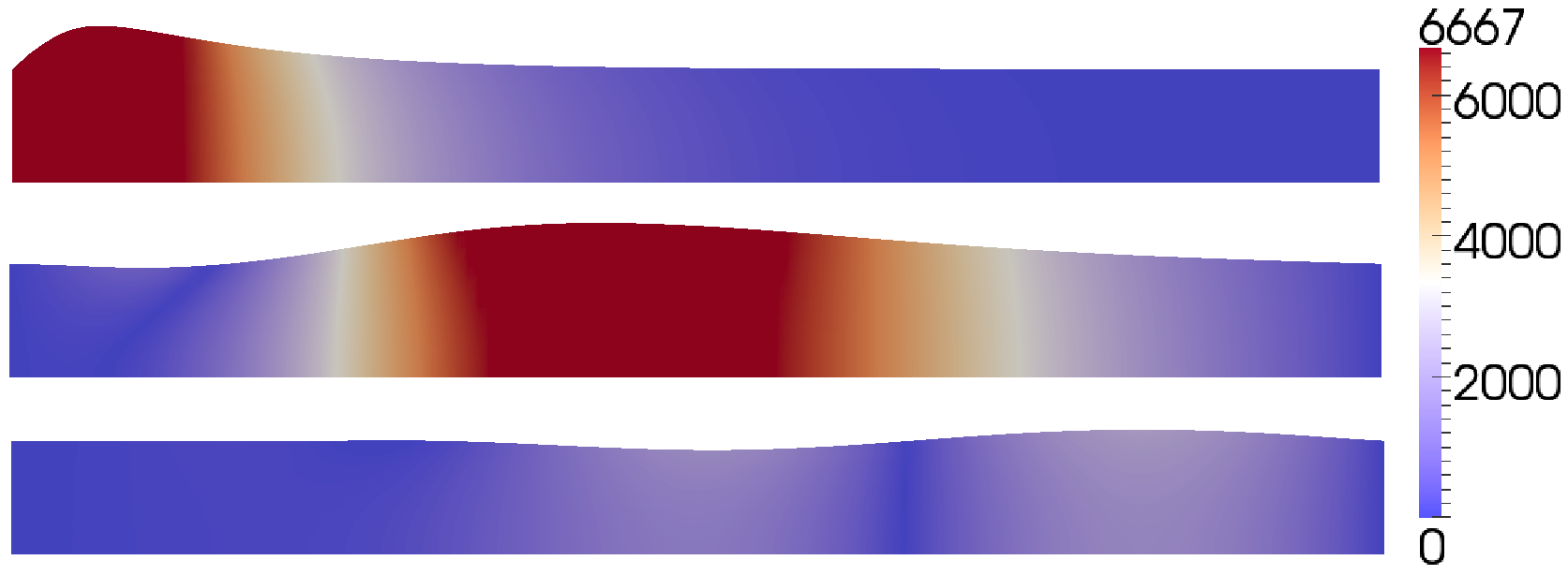}
\end{center}
\caption{Pressure surface plot on the deformed 2D fluid domain at times 1.5, 3.5 and 5.5 ms top to bottom. Units are $Pa$. The propagation of a wave is clearly visible. For visualization purposes, the deformation of the domain is magnified 100 times.}
\label{fig:0}
\end{figure}

\begin{figure}
\begin{center}
\includegraphics[width=0.3\textwidth]{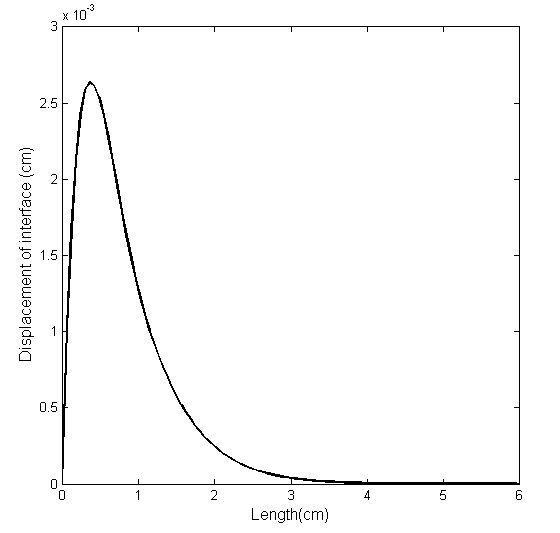} 
\includegraphics[width=0.3\textwidth]{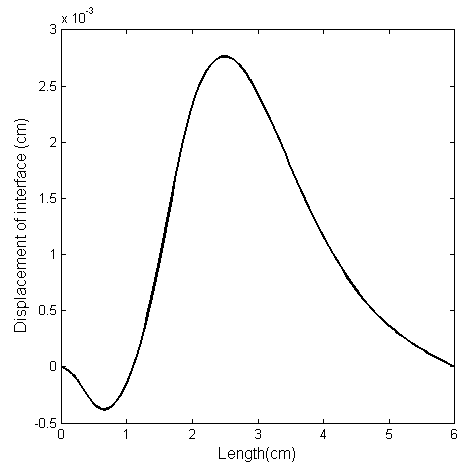} 
\includegraphics[width=0.3\textwidth]{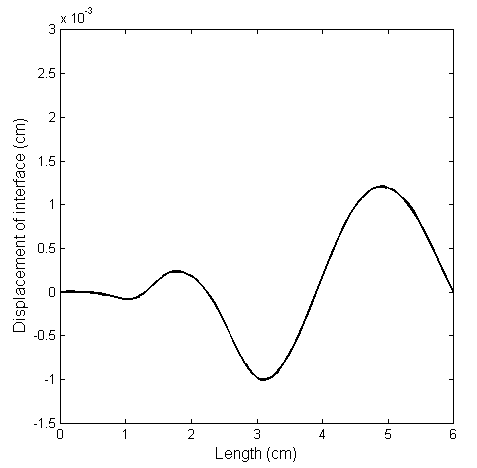}

\includegraphics[width=0.3\textwidth]{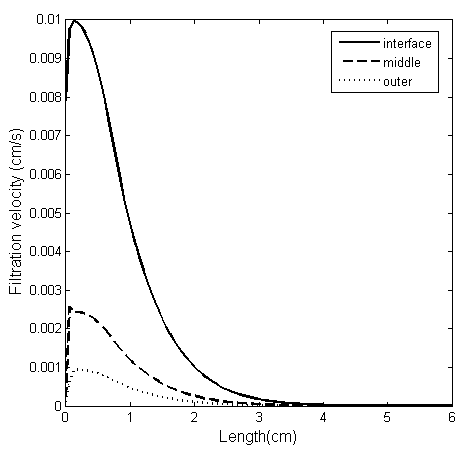} 
\includegraphics[width=0.3\textwidth]{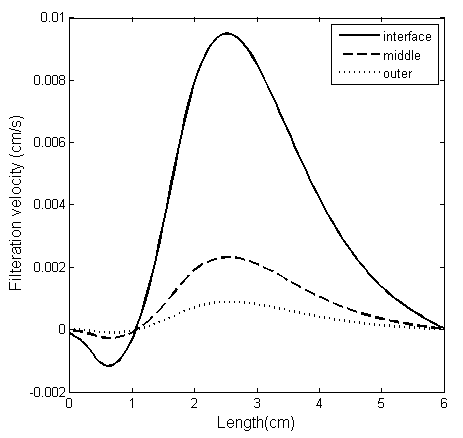} 
\includegraphics[width=0.3\textwidth]{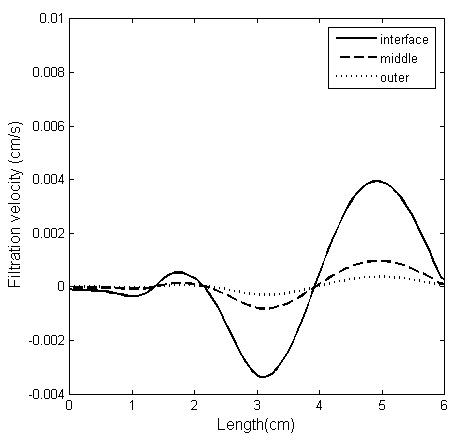}
\end{center}
\caption{Top panel: displacement of the fluid-structure interface at times 1.5, 3.5 and 5.5 ms from left to right. Bottom panel: intramural flow $\velp \cdot \nn$ at different planes in the arterial wall, located at the interface, at the intermediate section and at the outer layer.}
\label{fig:1}
\end{figure}

Some visualizations of the solution, calculated using the settings addressed below, are reported in Figures \ref{fig:0} and \ref{fig:1}. The former, qualitatively shows the propagation of a pressure wave along the channel, together with the corresponding deformation of the fluid domain at times $t_1=1.5,\,t_2=3.5,\,t_3=5.5$ ms. For visualization purposes, the vertical displacement is magnified 100 times in Figure \ref{fig:0}. In Figure \ref{fig:1}, top panel, we show the vertical displacement of the interface along the longitudinal axis of the channel. These plots show that the variable inflow pressure combined with the fluid-structure interaction mechanisms, generates a wave in the structure that propagates from left to right. This simulation is qualitatively similar to the ones obtained in \cite{bukac2013explicit} using a different discretization approach to model FSI, including an additional thin layer in the arterial wall. On the bottom panel, we show intramural flow $\velp \cdot \nn$ at different planes cutting the arterial wall in the longitudinal direction. These planes are located at the interface, at the intermediate section and at the outer layer. These plots show that the peaks of the intramural flow coincide with the ones of structure displacement and the corresponding peak of arterial pressure. Furthermore, we notice that the intramural velocity, $\velp \cdot \nn$, decreases as far as the fluid penetrates further into the wall. This is a consequence of equation \eqref{B3}, which prescribes that $\nabla \cdot \velp$ is not locally preserved, but depends on the rate of change in pressure and volumetric deformation of the structure. Indeed, this is how the poroelastic coupling shows up in the results.

\paragraph{Finite element approximation.} 
The proposed scheme is rather flexible with respect to the definition of the finite element spaces for each unknown, namely $\femvf,\,\fempf,\,\femvp,\,\fempp,\,\fems,\,\dotfems$ previously introduced. Owing to the weak enforcement of the interface conditions, in principle it does not require any restriction on the choice of the discrete approximation spaces on each side of the interface, provided that appropriate quadrature methods are adopted to calculate the integrals of test functions evaluated at the interface. Considering the heterogeneous nature of the equations that are coupled though the interface conditions, this is an interesting property of the scheme. Nevertheless, in the forthcoming numerical experiments we will consider some simple choices. To facilitate the exchange of information across the interface, we use the same finite elements for the approximation of the unknowns $\velf,\,\velp,\,\dispp$ that are combined in the interface conditions, see for example equation \eqref{CNF}. 

One of the drawbacks of using the Nitsche's method for interface conditions is that the stability of the scheme relies on the appropriate quantification of the penalty and stabilization operators, namely $\gamma_f,\,\gamma_{stab},\,\gamma_{stab}^{\prime}$. More precisely, these parameters must satisfy the requirements of Theorem \ref{th:stability} or Corollaries \ref{th:stability-1} and \ref{th:stability-2}, where the explicit expression of the constant $C_{TI}$ may not be available in general, but only their order of magnitude can be reasonably quantified. As a result, $\gamma_f,\,\gamma_{stab},\,\gamma_{stab}^{\prime}$ should be manually tuned for each specific application. For the particular case of FSI related to blood flow in arteries, the sensitivity of the numerical solution on $\gamma_f$ and $\gamma_{stab}$ has been carefully addressed in \cite{MR2498525} for a scheme similar to the one addressed here. According to these results combined with the insight provided by Corollary \ref{th:stability-2} on the choice of $\gamma_{stab}^{\prime}$, we have used the following values: $\gamma_f = 2500$, $\gamma_{stab} = 1$, $\gamma_{stab}^\prime = 0$. 

\paragraph{Convergence analysis.}
The objective of this section is to support the theoretical results on the accuracy of the proposed scheme, given in Theorem \ref{th:convergence}. The available theory shows that the main shortcomings of the scheme in terms of accuracy arise in the time approximation and are related to the splitting of the equations within each time step. For this reason, we limit our investigation to the variation of the time step $\tau$. As reference solution, since the analytical one is not available in this case, we use the numerical solution calculated using the monolithic scheme with a small time step equal to $\tau=10^{-6}$ s. This way, we make sure that the splitting error is not polluting the solution, and the approximation error related to the time discretization scheme is negligible. This solution will be denoted with the subscript $ref$. To guarantee a sufficient spatial resolution as well as \emph{inf-sup} stability, we use $\mathbb{P}^2-\mathbb{P}^1$ approximations for velocity and pressure in the blood flow, combined with $\mathbb{P}^2-\mathbb{P}^1$ approximation of Darcy's equation and $\mathbb{P}^2$ approximation of the structure displacement and velocity. We investigate the convergence properties of the scheme in the norm $\triple{\cdot}_{\heartsuit,\,N}^2$ that is used in Theorem \ref{th:convergence}. More precisely, we split it in four parts:
\begin{gather*}
\errorf^N :=  \rho_f\|\velf^N-\velf^{N,ref}\|^2_{L^2(\domf)},
\
\errors^N (a) :=  \rho_p\|\dotdispp^N - \dotdispp^{N,ref}\|^2_{L^2(\domp)},
\
\errors^N (c) := s_0\|\pp^N-\pp^{N,ref}\|^2_{L^2(\domp)},
\\
\errors^N (b) := 2 \mu_p \|{\boldsymbol D}(\dispp^N-\dispp^{N,ref})||^2_{L^2(\domp)}
+ \lambda_p ||\nabla \cdot (\dispp^N-\dispp^{N,ref})||^2_{L^2(\domp)},
\end{gather*}
corresponding to the fluid kinetic energy, the structure kinetic energy, the structure elastic stored energy and the pressure, respectively.
We calculate the error between the reference solution and solutions obtained using $\tau = \tau_0, \tau_0/2, \tau_0/4, \tau_0/8$ with $\tau_0 = 10^{-4}$ for simulations up to the final time $T=10^{-3}$ s. The mesh discretization step is $h=0.05$ cm for all cases.

\begin{table}
\begin{center}
\begin{tabular}{| l | c  c | c  c | c  c | c c |}
\hline
monolithic & $\sqrt{\errorf^N}$ & rate  & $\sqrt{\errors^N} (a)$ & rate & $\sqrt{\errors^N} (b)$  & rate & $\sqrt{\errors^N} (c)$ & rate \\
\hline
\hline
$\tau_0=10^{-4}$		&	2.14E-01	&		&	1.48E-01	&		&	5.24E-01	&		&	1.32E-02	&	\\\hline
$\tau_0/2$	&	1.05E-01	&	1.02	&	7.89E-02	&	0.91	&	2.82E-01	&	0.90	&	6.95E-03	&	0.92\\\hline
$\tau_0/4$	&	5.13E-02	&	1.04	&	4.03E-02	&	0.97	&	1.44E-01	&	0.97	&	3.53E-03	&	0.98\\\hline
$\tau_0/8$	&	2.45E-02	&	1.07	&	1.98E-02	&	1.03	&	7.07E-02	&	1.03	&	1.72E-03	&	1.03\\\hline
\end{tabular}

\medskip

\begin{tabular}{| l | c  c | c  c | c  c | c c |}
\hline
partitioned & $\sqrt{\errorf^N}$ & rate  & $\sqrt{\errors^N} (a)$ & rate & $\sqrt{\errors^N} (b)$  & rate & $\sqrt{\errors^N} (c)$ & rate \\
\hline
\hline															
$\tau_0=10^{-4}$	&	2.87E-01	&		&	1.84E-01	&		&	7.71E-01	&		&	1.96E-02	&	\\\hline
$\tau_0/2$	&	1.49E-01	&	0.94	&	9.91E-02	&	0.89	&	4.15E-01	&	0.90	&	1.01E-02	&	0.95\\\hline
$\tau_0/4$	&	7.58E-02	&	0.98	&	5.16E-02	&	0.94	&	2.13E-01	&	0.96	&	5.09E-03	&	0.99\\\hline
$\tau_0/8$	&	3.75E-02	&	1.01	&	2.59E-02	&	0.99	&	1.06E-01	&	1.01	&	2.49E-03	&	1.03\\\hline
\end{tabular}
\end{center}
\caption{Convergence in time of the monolithic and the partitioned scheme, respectively.}
\label{tab:2}
\end{table}

In Table \ref{tab:2} we show the convergence rate relative to the error indicators above calculated using both the monolithic and the loosely coupled scheme. We observe that, as expected, the error indicators scale as $C\tau$ when the monolithic scheme is used. This result confirms the statement of Assumption \ref{dip:consistency} with $p=1$, since we use the Backward Euler scheme for the time discretization of all the time dependent equations. Looking at the error of the loosely coupled scheme, we notice that for each of the indicators the magnitude of the error increases with respect to the monolithic scheme. This is the contribution of the splitting error. However, quite surprisingly, we observe that the total error of the loosely coupled scheme scales as $C\tau$, although Lemma \ref{lemma:spliterr} shows that the splitting error should asymptotically behave as $\sqrt{C_1 \frac \tau h + C_2 \tau^2 + C_3 \vartheta \tau}$. Combining these observations, we conclude that the splitting error is non negligible, but the dominant term among these three contributions is $\sqrt{C_2} \tau$. In other words, for the problem at hand, the available results suggest that $C_2 \gg C1,\,C_3$.

\paragraph{The loosely coupled scheme as a preconditioner for the monolithic formulation.}
In order to simplify the management of the algebraic degrees of freedom related to the finite element spaces, for this test case we adopt a $\mathbb{P}^1-\mathbb{P}^1$ approximations for velocity and pressure respectively. It is well known that this choice does not satisfy the \emph{inf-sup} stability condition \cite{MR2050138,MR1299729}. Resorting to a pressure stabilization method on the whole fluid domain is mandatory. Owing to its simplicity of implementation, we opt for the Brezzi-Pitkaranta scheme \cite{MR804083}, that is
\[
s_p(\pf^n,\bpf) := \gamma_p h^2 \int_{\domf} \nabla \pf^n \cdot \nabla \bpf d\boldsymbol x,
\]
where the stabilization parameter is selected as $\gamma_p =10^{-2}$ on the basis of numerical experiments.
The same type of spaces are used for the intramural filtration and pressure, namely $\velp,\pp$. 
Since equation \eqref{B3} is not enforcing the divergence-free constraint exactly, but the material turns out to be slightly compressible, equal order approximation is stable. We also use $\mathbb{P}^1$ finite elements to approximate the structure velocity and displacement.

Equation \eqref{eq:spectral} suggests that the partitioned scheme may act as an efficient preconditioner for the monolithic equation. On one hand, the corresponding matrices share similar spectral properties. On the other hand, the computational cost of solving the partitioned scheme is significantly lower than the one related to the monolithic formulation. Indeed, matrix $\maase$ features a block upper triangular structure that enables us to solve the structure, Darcy and blood flow problems independently. The performance of matrix $\maase$ used as a preconditioner of $\maasis$ is quantified by the numerical experiments reported in Table \ref{tab:3}. As an indicator of the system conditioning, we look at the number of GMRES iterations required to reduce below a given tolerance the relative residual, defined as the left hand side of \eqref{eq:convrate}. In the special case of positive definite matrices, a manipulation of \eqref{eq:convrate} allows us to show that the number of iterations (\# GMRES) required to reduce the relative residual of a factor $10^p$, can be estimated as \# GMRES $\simeq p \sqrt{\big(K(P^{-1}A)\big)}$, where $K(\cdot)$ is the spectral condition number. Since the initial relative residual is one, by definition, (\# GMRES) is equivalent to the number of iterations performed until the relative residual is less than $10^{-p}$. In the experiments that follow, we have used $p=6$. As a result, knowing that the conditioning of the FEM stiffness matrices scales as the square of the number of degrees of freedom, we expect that \# GMRES linearly scales with the number of degrees of freedom in absence of preconditioners. \emph{Optimal} preconditioners are those where the number of GMRES iterations becomes \emph{independent} of the dimension of the discrete problem.

\begin{table}
\begin{center}
Artery test case ($i$)

\smallskip 

\begin{tabular}{|c|c|c|c|}\hline
$\tau=10^{-4}$ & $h=0.05$ & $h=0.025$ & $h=0.0125$ \\\hline\hline
\# GMRES($\maasis$) & 211.4 & 446.4 & 1282.9 \\\hline
\# GMRES($(\maase)^{-1}\maasis$) & 10.9 & 12 & 13.9 \\\hline
\end{tabular}

\medskip

Artery test case ($i$)

\smallskip 

\begin{tabular}{|c|c|c|c|c|}\hline
$\tau=10^{-5}$ & $h=0.05$ & $h=0.025$ & $h=0.0125$ \\\hline\hline
\# GMRES($\maasis$) & 362.1 & 498.3 & 1194.4 \\\hline
\# GMRES($(\maase)^{-1}\maasis$) & 8 & 10 & 12.9 \\\hline
\end{tabular}

\medskip

Fractured reservoir test case ($ii$)

\smallskip 

\begin{tabular}{|c|c|c|c|c|}\hline
$\tau=10^{-3}$ & $h=1.06$ & $h=0.450$ & $h=0.212$ \\\hline\hline
\# GMRES($\maasis$) & 53.8 & 101.7 & 245.4 \\\hline
\# GMRES($(\maase)^{-1}\maasis$) & 2 & 2 & 2 \\\hline
\end{tabular}
\end{center}
\caption{Average number of GMRES iterations (\# GMRES) required to reduce the relative residual of a factor $10^{-6}$.
The values are calculated on the basis of the first 10 time steps of the simulation.}
\label{tab:3}
\end{table}

The results of Table \ref{tab:3} nicely agree with the general GMRES convergence theory, namely Theorem \ref{t:gmres}, and confirm that $\maase$ behaves as an optimal preconditioner for $\maasis$. Not only the number of iterations to solve the preconditioned system is nearly insensitive with respect to the mesh characteristic size, and consequently the number of degrees of freedom of the discrete problem, but the number of iterations is significantly smaller than in the non preconditioned case. Reminding that the inversion of $\maase$ is a relatively inexpensive operation, the preconditioned algorithm turns out to be a very effective solution method. Table \ref{tab:3} also suggests that the conditioning of the monolithic problem slightly increases when the time step is refined, especially for coarse meshes, while the good preconditioner performance seems to be unaffected. It indeed slightly improves, according to the fact that the loosely coupled scheme, and the related preconditioner, becomes more accurate and effective when the time discretization step decreases. We have tested this algorithm also using quadratic finite elements, $\mathbb{P}^2$, for all velocities and displacement fields. In the case of the coarsest mesh $h=0.05$ cm, GMRES converges in 613.5 (average) iterations, while solving the preconditioned only requires 11 iterations. The preconditioner seems to scale well also with respect to the FEM polynomial degree. The good results on preconditioner performance also correspond to a decrease in the computational time. For $h=0.05,\,\tau=10^{-4}$, the calculation of $60$ time steps of the monolithic scheme require 4.73 seconds, while for the preconditioned method the time is 1.85 seconds. For $\tau=10^{-4}$ and $600$ time steps the computational times are respectively, 65.8 seconds and 14.6 seconds. Finally, we observe that these considerations apply also to the test case ($ii$), addressed in the next section, which represents a prototype problem for flow though a fractured reservoir.

%>>>>>>>>>>>>>>>>>>>>>>>>>>>>>>>>>>>>>>>>>>>>>>>>>>>>>>>>>>>>>>>>>>>>>>>>>>>>>>>>>>>>>>>>>>>>>>>>>>>>>>>>>>>>>>>>>>>>>>>>>>>>>>>>>

\subsection{Application of the scheme to flow through fractured reservoirs}

\begin{figure}
\begin{center}
\includegraphics[width=\textwidth]{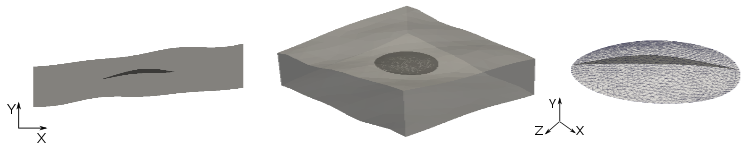}
\end{center}
\caption{Left: Physical computational domain $\Omega=\Omega^f \cup \Omega^p$ in 2D.  Middle: Physical domain $\Omega^p$ in 3D. Right: Physical domain $\Omega^f$ in 3D.  }
\label{FigDom}
\end{figure}

In this example we focus on fluid-structure interaction in the context of modelling the interaction between a stationary fracture filled with fluid and the surrounding poroelastic medium. We consider a general case where the hydraulic conductivity is a tensor given by $\boldsymbol \kappa = \frac{\boldsymbol K}{\mu_f}$, where $\boldsymbol K$ is the permeability tensor. Since the  effect of the porous media on the conduit flow at the interface is more prominent in this application, we replace the no-slip condition \eqref{CBJS} by the Beavers-Joseph-Saffman condition \eqref{CBJS_bis} with the coefficient 
\begin{equation}
\beta = \left(\frac{\alpha \mu_f \sqrt{3}}{\sqrt{\textrm{tr} (\boldsymbol{K})}} \right)^{-1}.
\end{equation}
We consider 2D and 3D test cases. In 2D, our reference domain is a square $[-100,100]^2$ m. A fracture is positioned in the middle of the square, whose boundary is given by $$\hat{y}^2 = 0.008^2(\hat{x}-35)^2(\hat{x}+35^2).$$
The fracture represents the reference fluid domain $\hat{\Omega}^f$, while the reference poroelastic structure domain is given as $\hat{\Omega}^p=\hat{\Omega} \backslash \hat{\Omega}^f.$
Following the approach in~\cite{ganis1}, to obtain a more realistic domain, we transform the reference domain $\hat{\Omega}$ onto the physical domain $\Omega$ via the mapping
\begin{equation}
\left[\begin{array}{c}
x \\ y 
\end{array} \right] (\hat{x},\hat{y}) = 
\left[\begin{array}{c}
\hat{x} \\
5\cos\left(\displaystyle\frac{\hat{x}+\hat{y}}{100}\right)\cos\left(\displaystyle\frac{\pi\hat{x}+ \hat{y}}{100}\right)^2+\hat{y}/5-\hat{x}/10
\end{array} \right].
\end{equation}
The extension to the 3D test case consists of introducing a reference domain  $\hat{\Omega} = \hat{\Omega}^f \cup \hat{\Omega}^p$ which is a cube $[-100,100]^3$ m. A circular fracture with radius of $35$ m and maximum height of $20$ m is positioned in the middle of the cube. 
Then, we apply the mapping
\begin{equation}
\left[\begin{array}{c}
x \\ y \\ z 
\end{array} \right] (\hat{x},\hat{y},\hat{z}) = 
\left[\begin{array}{c}
\hat{x} \\
5\cos\left(\displaystyle\frac{\pi\hat{z}+\hat{x}+\hat{y}}{100}\right)\cos\left(\displaystyle\frac{\pi\hat{x}+ \hat{y}}{100}\right)^2+\hat{y}/5+\hat{x}(\hat{z}/100-1)/10 \\  \hat{z} 
\end{array} \right].
\end{equation}
Figure~\ref{FigDom} shows the physical domain $\Omega$ in 2D (left) and 3D (right).  

The flow is driven by the injection of the fluid into the fracture with the constant rate $g = 25$ kg/s. The fluid is injected into a circle (2D) or a cube (3D) of radius $7$ m in the center of the reference fracture $\hat{\Omega}^f$. On all external boundaries, we prescribe the no flow  condition $\boldsymbol q \cdot \boldsymbol n=0$, zero normal displacement $\boldsymbol U \cdot \boldsymbol n=0$, and zero shear traction $\ttau \cdot \boldsymbol \sigma_f \nn =0 $. The computational time is $T=5.1$ hours. The problem was solved using the time step $\Delta t = 0.1$ seconds. The remaining parameters are given in Table~\ref{Tgeo}. These parameters are similar to the ones in~\cite{ganis1}, and are taken to be within realistic range for applications to hydraulic fracturing.   
\begin{table}[ht!]
\begin{center}
{\small{
\begin{tabular}{l c c c}
\hline
\textbf{Parameter} & \textbf{Symbol} & \textbf{Units} &\textbf{Value} \\
\hline
\hline
\textbf{Poroelastic wall density}& $\rho_{p}$ & (kg/m$^3$) & $897$  \\
\textbf{Fluid density}& $\rho_f$ & (kg/m$^3$)& $897$\\
\textbf{Dyn. viscosity}& $\mu$ & (Pa s) & $10^{-3}$    \\
\textbf{Lam\'e coeff.}& $\mu_p $ & (Pa) & $2.92 \times 10^{10}$  \\
\textbf{Lam\'e coeff.}& $\lambda_p $ & (Pa) & $1.94 \times 10^{10}$ \\
\textbf{Hydraulic conductivity (2D)}& $\boldsymbol \kappa $ & (m$^2$/Pa s) & $diag(200,50) \times 10^{-12}$ \\
\textbf{Hydraulic conductivity (3D)}& $\boldsymbol \kappa $ & (m$^2$/Pa s) & $diag(200,50,200) \times 10^{-12}$ \\
\textbf{Mass storativity coeff.}& $s_0 $ & (Pa$^{-1}$) & $6.9 \times 10^{-5}$   \\
\textbf{Biot-Willis constant}& $\alpha$ &  & 1 \\
\textbf{Beavers-Joseph-Saffman coefficient (2D)}& $\beta $ & (m$^2$/Pa s) & $ 2.88\cdot 10^{-4}$  \\
\textbf{Beavers-Joseph-Saffman coefficient (3D)}& $\beta $ & (m$^2$/Pa s) & $ 3.88\cdot 10^{-4}$  \\
\textbf{Total simulation time}& $T$ & (s) & $ 18 600$  
 \\\hline
\end{tabular}
}}
\end{center}
\caption{Fluid and structure parameters.}
\label{Tgeo}
\end{table}

In the 2D case, we use the conforming meshes for the discretization of the fluid and structure domains. In 3D case, we exploit the weak enforcement of the interface conditions and use non-conforming computational meshes. We use the preconditioned monolithic numerical solver. Its implementation in \textit{FreeFem++} facilitates the interpolation of finite element functions across the interface. We adopt the $\mathbb{P}^1-\mathbb{P}^1$ approximation for the fluid velocity and pressure, complemented with the pressure stabilization $s_p(p_{f,h}, \psi_{f,h})$. For the Darcy fluid-pressure approximation in the poroelastic medium, as for the structure velocity and displacement,  we again use the $\mathbb{P}^1-\mathbb{P}^1$ finite elements. However, due to the large time-step in this example, we are numerically closer to the divergence-free regime. Thus, we add a Darcy pressure stabilization given by
$$
s_q(p_{p,h}, \psi_{p,h}) = \gamma_{q} h^2 \int_{\Omega_p} \nabla p_{p,h} \cdot \nabla \psi_{p,h} d \boldsymbol {x},
$$
where the stabilization parameter is selected as $\gamma_q = 10^{-3}$. 
Finally, the choice of the penalty and stabilization parameters is  $\gamma_f = 1500$, $\gamma_{stab} = 1$, $\gamma_{stab}^\prime = 0$. 

\begin{figure}
\begin{center}
\includegraphics[width=1\textwidth]{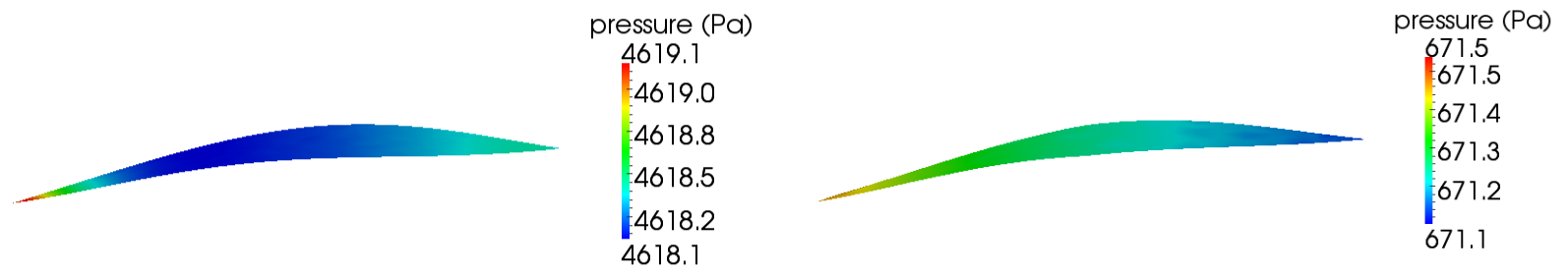}
\includegraphics[width=1\textwidth]{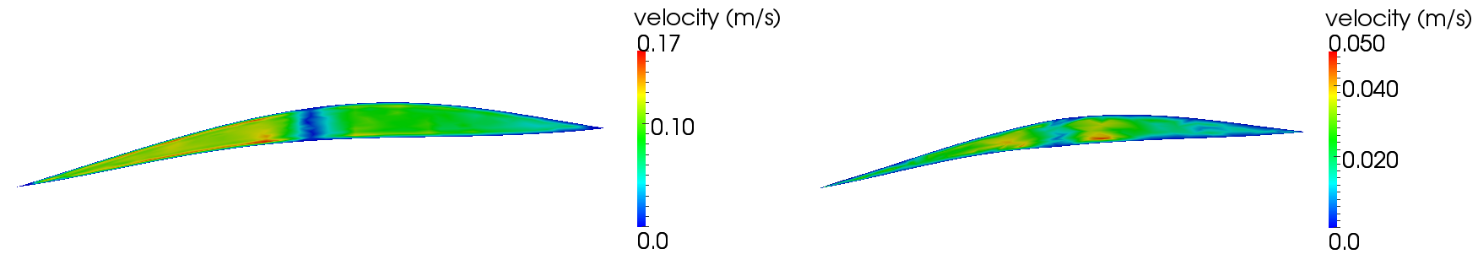}
\end{center}
\caption{Top: Fluid pressure in 2D (left) and a slice of the pressure in the middle of the fracture in 3D (right)  at final time $T$. 
Bottom: Fluid velocity in 2D (left) and  and a slice of the velocity in the middle of the fracture in 3D (right)  at final time $T$.}
\label{FigFracture23D}
\end{figure}

\begin{figure}
\begin{center}
\includegraphics[width=0.45\textwidth]{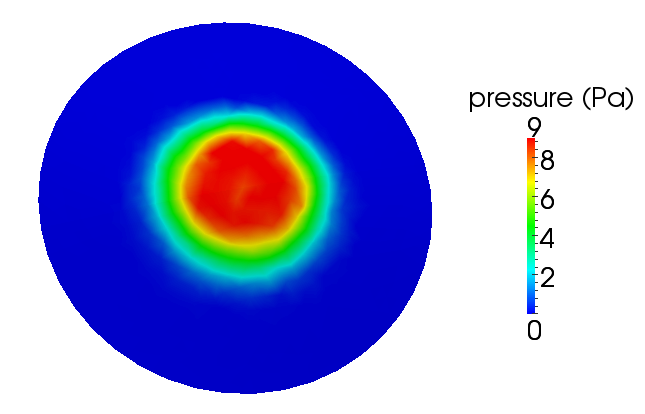}
\includegraphics[width=0.45\textwidth]{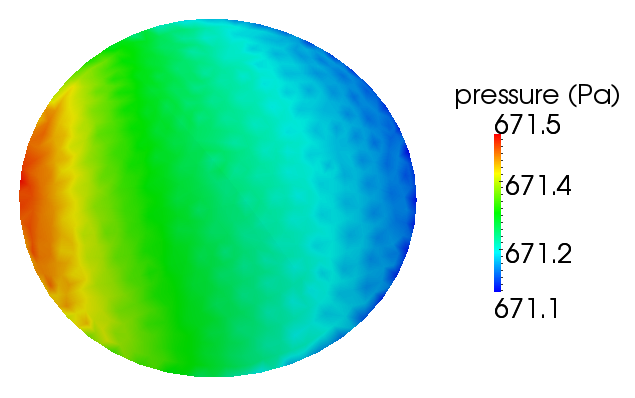}
\end{center}
\caption{Fluid pressure in 3D at time $t=5$ seconds (left) and at the final time $T$ (right). }
\label{FigFracture3D}
\end{figure}

\begin{figure}
\begin{center}
\includegraphics[width=1\textwidth]{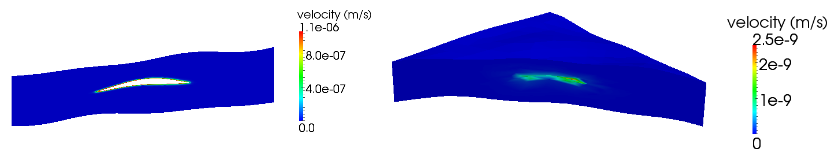}
\end{center}
\caption{Fluid velocity in the surrounding medium in 2D  (left) and 3D (right) at final time $T$. }
\label{FigDarcyVelocity}
\end{figure}

\begin{figure}
\begin{center}
\includegraphics[width=0.45\textwidth]{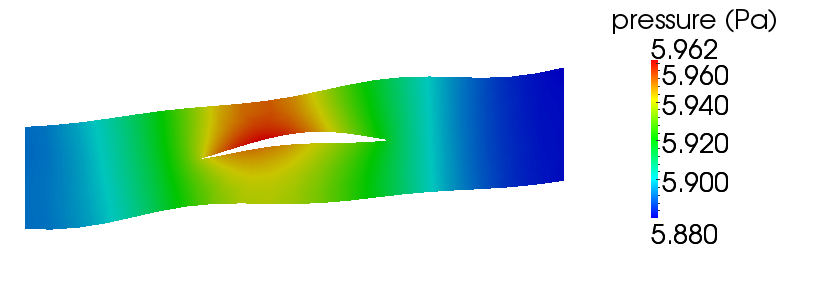}
\includegraphics[width=0.45\textwidth]{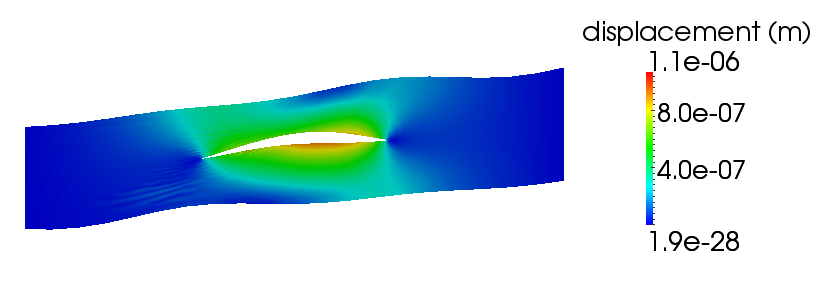}
\end{center}
\caption{Fluid pressure (left) and displacement (right) in surrounding medium in 2D at final time $T$. }
\label{FigFracture2D}
\end{figure}

\begin{figure}
\begin{center}
\includegraphics[width=\textwidth]{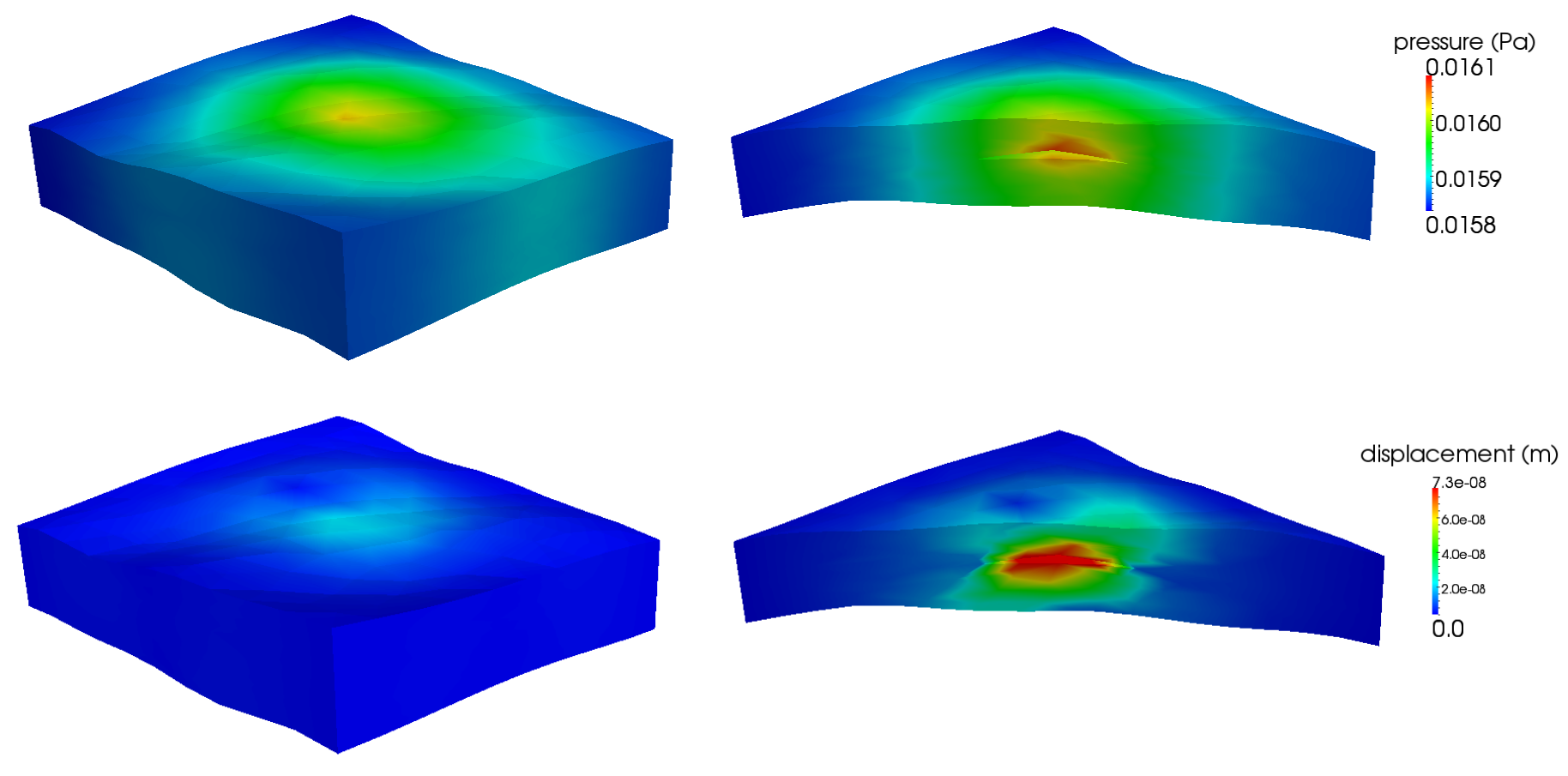}
\end{center}
\caption{Top: Pressure in the poroelastic medium at final time $T$. Bottom: Displacement in the poroelastic medium at final time $T$. }
\label{FigRock}
\end{figure}

\begin{figure}
\begin{center}
\includegraphics[width=1.05\textwidth]{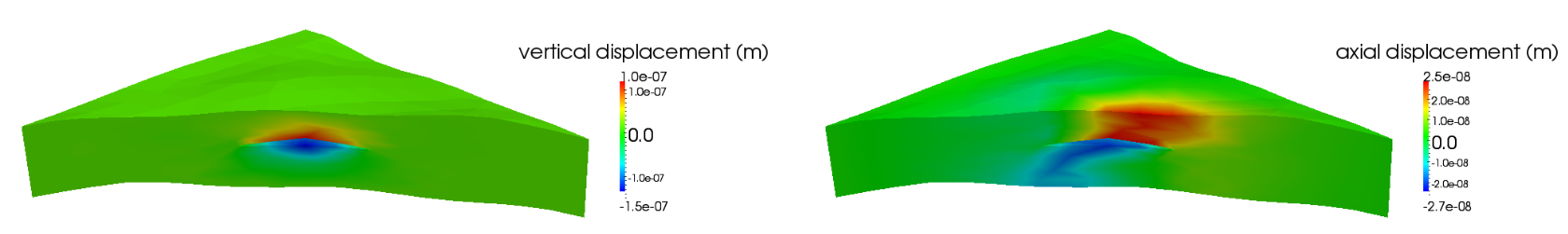}
\end{center}
\caption{Displacement of the poroelastic medium in 3D in the vertical  (left) and axial (right) direction at final time $T$. }
\label{FigRockDispl}
\end{figure}

The visualization of the injection of the fluid into the fracture, and the fluid leakoff to the surrounding rock at the final time $T$ is shown in a sequence of Figures, numbered from \ref{FigFracture23D} to \ref{FigRockDispl}.

Figure~\ref{FigFracture23D} shows the pressure (top) and the velocity (bottom) in the fracture for the 2D (left) and 3D (right) test cases, at final time $T$. The 3D plots are obtained by taking a slice through the center of the domain. A top view of the fluid pressure in the fracture is shown in Figure~\ref{FigFracture3D} (as the pressure range is very narrow in the right panel, small mesh dependent fluctuations of the pressure are visible from the pressure contour plot). At the beginning of the process the simulations capture the expected local pressure increase in the region of fluid injection, while at the final time, the fluid pressure is the largest at the tip of the crack in both 2D and 3D test cases. Figure~\ref{FigDarcyVelocity} shows the comparison of Darcy velocity in the surrounding medium obtained in 2D (left) and 3D (right)  at final time $T$. Pressure and displacement are shown in the surrounding medium, in the 2D case in Figure~\ref{FigFracture2D}, and in the 3D test case  in Figure~\ref{FigRock}. Finally, the displacement of the surrounding medium in the vertical (left) and axial (right) direction is shown in Figure~\ref{FigRockDispl}.

This example demonstrates the ability of our algorithm to handle complex, three-dimensional simulations in different applications.  
Retaining the advantages of a monolithic scheme, the computational time is significantly reduced due to a well-designed preconditioner. Our model includes Darcy equations in the mixed formulation which are necessary to compute accurately the Darcy velocity in the surrounding rock. Furthermore, we are able to capture displacement of the fracture which has a complex, very slender domain with sharp edges. 

%>>>>>>>>>>>>>>>>>>>>>>>>>>>>>>>>>>>>>>>>>>>>>>>>>>>>>>>>>>>>>>>>>>>>>>>>>>>>>>>>>>>>>>>>>>>>>>>>>>>>>>>>>>>>>>>>>>>>>>>>>>>>>>>>>

\section{Conclusions}

We have studied the interaction of a free fluid with a poroelastic material. After setting appropriate governing equations on adjacent domains and discussing the corresponding interface conditions, we have considered the discretization of the problem in the framework of the finite element method. Particular attention must be devoted to the approximation of the interface conditions, which are non standard with respect to the ones that arise in the coupling of homogenous partial differential equations. We have shown that the Nitsche's method, well known the weak approximation of boundary and interface conditions for elliptic or parabolic problems, is very appropriate to enforce the interface constraints in the variational formulation. The resulting discrete, fully coupled problem features good stability properties. Since all the interface conditions correspond to suitable operators in the variational problem, time-lagging allows to split the fully coupled problem into subproblems, relative to the main governing equations, such as free fluid flow, Darcy filtration and elastodynamics. Also the resulting loosely coupled problem formulation turns out to be stable, provided that it is combined with suitable stabilization operators. This solution approach is very effective from the computational standpoint, but suffers from low accuracy, as we have pointed out in the error analysis of the scheme. In order to merge the computational efficiency of the loosely coupled scheme with the good accuracy and stability properties of the monolithic formulation, we develop a numerical solver where the former scheme acts as a preconditioner for the latter. The theory and the numerical results suggest that this approach is very effective because the loosely coupled scheme behaves as an optimal preconditioner for the monolithic formulation. This solution algorithm turns out to be very robust with respect to the characteristic physical parameters of the problem. Indeed, we have successfully applied it to the analysis of a problem related to blood flow in arteries as well as to the study of subsurface flow and deformation of a fractured reservoir.

%>>>>>>>>>>>>>>>>>>>>>>>>>>>>>>>>>>>>>>>>>>>>>>>>>>>>>>>>>>>>>>>>>>>>>>>>>>>>>>>>>>>>>>>>>>>>>>>>>>>>>>>>>>>>>>>>>>>>>>>>>>>>>>>>>

\section*{References}
%\bibliographystyle{elsarticle-num}
%\bibliographystyle{plain}
%\bibliography{fsi-biot-nitsche-new}

%>>>>>>>>>>>>>>>>>>>>>>>>>>>>>>>>>>>>>>>>>>>>>>>>>>>>>>>>>>>>>>>>>>>>>>>>>>>>>>>>>>>>>>>>>>>>>>>>>>>>>>>>>>>>>>>>>>>>>>>>>>>>>>>>>

\end{document}